\documentclass[12pt]{amsart}
\usepackage{amsmath, amssymb, latexsym, amsthm, mathrsfs, geometry}
\usepackage[all]{xy}
\newcommand{\ed}{

\subsection*{Acknowledgments}
We thank Daniel Hathaway for his comments on the proof of Theorem \ref{thm:cflr}.
A part of this material was presented at the Tel Aviv University
Set Theory Seminar. We thank Moti Gitik and Assaf Rinot for their useful feedback.
The research of the first named author was partially supported by the Science
Absorption Center, Ministry of Immigrant Absorption, the State of Israel.
The research of the second named author was partially supported by
the Israel Science Foundation, grant number 1053/11.
This is publication 1032 of the second named author.

\end{document}}

      \newenvironment{changemargin}[2]{\begin{list}{}{
         \setlength{\topsep}{0pt}\setlength{\leftmargin}{0pt}
         \setlength{\rightmargin}{0pt}
         \setlength{\listparindent}{\parindent}
         \setlength{\itemindent}{\parindent}
         \setlength{\parsep}{0pt plus 1pt}
         \addtolength{\leftmargin}{#1}\addtolength{\rightmargin}{#2}
         }\item }{\end{list}}

\newcommand{\nc}{\newcommand}
\newcommand{\two}{\{0,1\}}
\nc{\blue}[1]{{\color{blue} #1}}
\nc{\productive}[2]{\bigl(#1,\allowbreak #2\bigr)^\x}
\nc{\name}[1]{\dot{#1}}
\nc{\forces}{\Vdash}
\nc{\set}[2]{\left\{\,#1 : #2\,\right\}}
\nc{\smallset}[2]{\bigl\{\,#1 : #2\,\bigr\}}
\nc{\seq}[2]{\left\la\, #1 : #2\,\right\ra}
\nc{\smallseq}[2]{\la\, #1 : #2\,\ra}
\nc{\cube}{(\Cantor)^\N}
\nc{\dP}{\dot{\bbP}}
\nc{\dQ}{\dot{\bbQ}}
\nc{\Match}{\op{Match}}
\nc{\concat}[1]{\hat{\phantom{a}}\langle #1\rangle}
\nc{\poset}{\mathbb{P}}
\nc{\fn}[1]{{\op{Fn}(#1\times\w,2)}}
\nc{\linadd}{\op{linadd}}
\nc{\nonprod}{\non^\x}
\nc{\Ga}{\Gamma}
\nc{\Om}{\Omega}
\nc{\alephes}{{\aleph_0}}
\nc{\my}[1]{{\color{red} #1}}
\nc{\Cp}{\op{C}_p}
\nc{\Bp}{\op{B}_p}
\nc{\Pa}[8]{\bibitem{#1} {#2}, \emph{#3}, {#4} \textbf{#5} ({#6}), {#7}--{#8}.}
\nc{\tPa}[5]{\bibitem{#1} {#2}, \emph{#3}, {#4}, to appear.}
\nc{\sPa}[4]{\bibitem{#1} {#2}, \emph{#3}, {#4}, submitted.}
\nc{\Bc}[9]{\bibitem{#1} {#2}, \emph{#3}, in: \textbf{#4} (#5), #6 #7, #8--#9.}
\nc{\fD}{\mathfrak{D}}
\nc{\fX}{\mathfrak{X}}
\nc{\Onbd}{\Op_{\mathrm{nbd}}} 
\nc{\Omnb}{\Om_{\mathrm{nbd}}} 
\nc{\od}{\mathfrak{od}}
\nc{\Setting}[7]{\xymatrix@R=4pt@C=7pt{#1\ar@{-}[r]&#2\ar@{-}[r]&#3\\&#4\ar@{-}[u]\\
#5\ar@{-}[uu]\ar@{-}[r] & #6\ar@{-}[u]\ar@{-}[r] & #7\ar@{-}[uu]}}
\nc{\mx}[1]{\begin{matrix}#1\end{matrix}}
\nc{\plim}{p\txt{-}\lim}
\nc{\Bgp}{{\Z^\N}}
\nc{\Cgp}{{{\Z_2}^\N}}
\nc{\Cite}[1]{\textbf{[#1]}}
\nc{\Next}[1]{{#1^+}}
\nc{\Fr}{\mathit{F\!r}}
\nc{\intvl}[2]{{\bigl[#1(#2),\allowbreak #1(#2\!+\!1)\bigr)}}
\nc{\bintvl}[2]{{\bigl[#1, #2\bigr)}}
\nc{\Bdd}{\mathbf{B}}
\nc{\Ax}{\mathsf{Ax}}
\nc{\Dfin}{\mathfrak{D}_\mathrm{fin}}
\nc{\grbl}{{\mbox{\textit{\tiny gp}}}}
\nc{\bbP}{\mathbb{P}}
\nc{\bbC}{\mathbb{C}}
\nc{\bbD}{\mathbb{D}}
\nc{\bbX}{\mathbb{X}}
\nc{\bbO}{\mathbb{O}}
\nc{\bbQ}{\mathbb{Q}}
\nc{\BOfat}{\B_{\Om_{\mathrm{fat}}}}
\nc{\Bgood}{\B_{\mathrm{good}}}
\nc{\compactN}{\cl{\mathbb{N}}}
\nc{\blocks}[2]{\op{cl}_{#2}(#1)}
\nc{\blocksplus}[2]{\op{cl}^+_{#2}(#1)}
\nc{\arx}[1]{\texttt{http://arxiv.org/math/#1}}
\nc{\bq}{\begin{quote}}
\nc{\eq}{\end{quote}}
\nc{\cl}[1]{\overline{#1}}
\nc{\CH}{the Continuum Hypothesis}
\nc{\MA}{Martin's Axiom}
\nc{\Bfat}{\B_\mathrm{fat}}
\nc{\inv}{^{-1}}
\nc{\Cantor}{{\two^\N}}
\nc{\bP}{\mathbf{P}}
\nc{\Fn}{\op{Fn}}
\nc{\bof}{\op{\fb}}
\nc{\bofF}{\bof(\cF)}
\nc{\sr}[3]{{\txt{#1\\#3}}}
\nc{\smallbinom}[2]{\bigl(\begin{smallmatrix}#1\\#2\end{smallmatrix}\bigr)}
\nc{\gp}{\binom{\Om}{\Ga}}
\nc{\gpsmall}{\mbox{$\gp$}}
\nc{\gig}{\gimel}
\nc{\gns}{\sone(\Om,\gig)}
\nc{\nsr}[2]{#1}
\nc{\N}{\w}
\nc{\NN}{{\N^{\N}}}
\nc{\ZN}{{\Z^{\N}}}
\nc{\NNup}{{\N^{\uparrow\N}}}
\nc{\PN}{{P(\N)}}
\nc{\roth}{{[\w]^{\w}}} 
\nc{\Fin}{[\N]^{\mbox{\tiny $<\!\infty$}}} 
\nc{\ici}{{[\w]^{(\w,\w)}}}
\nc{\Inc}{{\compactN^{\uparrow\N}}}
\nc{\powInc}[1]{{\big(\Inc\big)^{#1}}}
\nc{\powFin}[1]{{\big(\Fin\big)^{#1}}}
\nc{\powPN}[1]{{\big(\PN\big)^{#1}}}
\nc{\NcompactN}{{\compactN^\N}}
\nc{\setseq}[1]{\{#1 : n\in\N\}}
\nc{\sseq}[1]{\{#1 : n\in\N\}}
\nc{\Uarrow}{\smash{\big\uparrow}}
\nc{\LE}{\preccurlyeq}
\nc{\GE}{\succcurlyeq}
\nc{\op}{\operatorname}
\nc{\im}{\op{im}}
\nc{\Span}{\op{span}}
\nc{\maxfin}{\op{maxfin}}
\nc{\ran}{\op{range}}
\nc{\iso}{\cong}
\nc{\Madd}{{\M}^\star}
\nc{\cI}{\mathcal{I}}
\nc{\cJ}{\mathcal{J}}
\nc{\scrA}{\mathscr{A}}
\nc{\scrB}{\mathscr{B}}
\nc{\scrC}{\mathscr{C}}
\nc{\scrD}{\mathscr{D}}
\nc{\scrF}{\mathscr{F}}
\nc{\scrK}{\mathscr{K}}
\nc{\A}{\forall}
\nc{\B}{\mathrm{B}}
\nc{\cB}{\mathcal{B}}
\nc{\cR}{\mathcal{R}}
\nc{\bB}{\mathbf{B}}
\nc{\BG}{\B_\Ga}
\nc{\BL}{\B_\Lambda}
\nc{\BT}{\B_\Tau}
\nc{\BTstar}{\B_{\Tau^*}}
\nc{\BO}{\B_\Om}
\nc{\DO}{\cD_\Om}
\nc{\KO}{\cK_\Om}
\nc{\CG}{C_\Ga}
\nc{\CL}{C_\Lambda}
\nc{\CT}{C_\Tau}
\nc{\CTstar}{C_{\Tau^*}}
\nc{\CO}{C_\Om}
\nc{\COgp}{C_{\Om^{\grbl}}}
\nc{\CLgp}{C_{\Lambda^{\grbl}}}
\nc{\BOgp}{\B_{\Om}^{\grbl}}
\nc{\BLgp}{\B_{\Lambda^{\grbl}}}
\nc{\sfC}{\mathsf{C}}
\nc{\sfD}{\mathsf{D}}
\nc{\bD}{\mathbf{D}}
\nc{\Tau}{\mathrm{T}}
\nc{\cA}{\mathcal{A}}
\nc{\cC}{\mathcal{C}}
\nc{\cK}{\mathcal{K}}
\nc{\cD}{\mathcal{D}}
\nc{\cF}{\mathcal{F}}
\nc{\cS}{\mathcal{S}}
\nc{\cG}{\mathcal{G}}
\nc{\cY}{\mathcal{Y}}
\nc{\J}{\mathcal{J}}
\nc{\cL}{\mathcal{L}}
\nc{\cM}{\mathcal{M}}
\nc{\cN}{\mathcal{N}}
\nc{\cO}{\mathcal{O}}
\nc{\Op}{\mathrm{O}}
\nc{\rmA}{\mathrm{A}}
\nc{\rmB}{\mathrm{B}}
\nc{\cP}{\mathcal{P}}
\nc{\Q}{\mathbb{Q}}
\nc{\R}{\mathbb{R}}
\nc{\cU}{\mathcal{U}}
\nc{\Un}{\bigcup}
\nc{\cV}{\mathcal{V}}
\nc{\cW}{\mathcal{W}}
\nc{\Z}{{\mathbb Z}}
\nc{\Impl}{\Rightarrow}
\long\def\forget#1\forgotten{}
\nc{\ft}{\mathfrak{t}}
\nc{\fb}{\mathfrak{b}}
\nc{\fc}{\mathfrak{c}}
\nc{\fd}{\mathfrak{d}}
\nc{\fg}{\mathfrak{g}}
\nc{\oo}{\infty}
\nc{\fr}{\mathfrak{r}}
\nc{\fu}{\mathfrak{u}}
\nc{\fh}{\mathfrak{h}}
\nc{\fp}{\mathfrak{p}}
\nc{\fq}{\mathfrak{q}}
\nc{\lr}{\mathfrak{lr}}
\nc{\lx}{\mathfrak{lx}}
\nc{\fj}{\mathfrak{j}}
\nc{\fs}{\mathfrak{s}}
\nc{\w}{\omega}
\nc{\x}{\times}
\nc{\Iff}{\Leftrightarrow}

\nc{\nin}{\notin}
\nc{\cat}{\hat{\ }}
\nc{\sub}{\subseteq}
\nc{\supp}{\op{supp}}
\nc{\spst}{\supseteq}
\nc{\sm}{\setminus}
\nc{\as}{\subseteq^*}
\nc{\rest}{\restriction}

\nc{\la}{\langle}
\nc{\ra}{\rangle}
\nc{\dom}{\op{dom}}
\nc{\cov}{\op{cov}}
\nc{\add}{\op{add}}
\nc{\cof}{\op{cof}}
\nc{\cf}{\op{cf}}
\nc{\non}{\op{non}}
\nc{\unif}{\op{non}}
\nc{\COV}{\op{COV}}
\nc{\ADD}{\op{ADD}}
\nc{\COF}{\op{COF}}
\nc{\NON}{\op{NON}}
\nc{\impl}{\to}
\nc{\card}[1]{\left|#1\right|}
\nc{\Wlog}{without loss of generality}
\newtheorem{thm}{Theorem}[section]
\nc{\bthm}{\begin{thm}} \nc{\ethm}{\end{thm}}
\newtheorem{prop}[thm]{Proposition}
\nc{\bprp}{\begin{prop}} \nc{\eprp}{\end{prop}}
\newtheorem{fact}[thm]{Fact}
\nc{\bfct}{\begin{fact}} \nc{\efct}{\end{fact}}
\newtheorem{prob}[thm]{Problem}
\nc{\bprb}{\begin{prob}} \nc{\eprb}{\end{prob}}
\newtheorem{lem}[thm]{Lemma}
\nc{\blem}{\begin{lem}} \nc{\elem}{\end{lem}}
\newtheorem{claim}[thm]{Claim}
\nc{\bclm}{\begin{claim}} \nc{\eclm}{\end{claim}}
\newtheorem{cor}[thm]{Corollary}
\nc{\bcor}{\begin{cor}} \nc{\ecor}{\end{cor}}
\newtheorem{conj}[thm]{Conjecture}
\nc{\bcnj}{\begin{conj}} \nc{\ecnj}{\end{conj}}
\theoremstyle{definition}
\newtheorem{defn}[thm]{Definition}
\nc{\bdfn}{\begin{defn}} \nc{\edfn}{\end{defn}}
\theoremstyle{remark}
\newtheorem{rem}[thm]{Remark}
\nc{\brem}{\begin{rem}} \nc{\erem}{\end{rem}}
\newtheorem{cnv}[thm]{Convention}
\nc{\bcnv}{\begin{cnv}} \nc{\ecnv}{\end{cnv}}
\newtheorem{exam}[thm]{Example}
\nc{\bexm}{\begin{exam}} \nc{\eexm}{\end{exam}}
\nc{\bpf}{\begin{proof}} \nc{\epf}{\end{proof}}
\nc{\be}{\begin{enumerate}}
\nc{\ee}{\end{enumerate}}
\nc{\bi}{\begin{itemize}}
\nc{\itm}{\item}
\nc{\ei}{\end{itemize}}
\nc{\Subsection}[1]{\goodbreak\subsection*{#1}}
\nc{\ev}[1]{\left[ #1 \right]}
\nc{\beq}{\begin{eqnarray*}}\nc{\eeq}{\end{eqnarray*}}

\forget
\forgotten

\nc{\sone}{\mathsf{S}_1}
\nc{\sfin}{\mathsf{S}_\mathrm{fin}}
\nc{\ufin}{\mathsf{U}_\mathrm{fin}}
\nc{\Split}{\mathsf{Split}}

\nc{\gone}{\mathsf{G}_1}    \nc{\gfin}{\mathsf{G}_\mathrm{fin}}

\title[Linear refinements and selections]{The linear refinement number and selection theory}

\author{Micha\l{} Machura}
\address[Machura]{Institute of Mathematics, University of Silesia, ul.\ Bankowa 14, 40-007 Katowice, Poland; and
Department of Mathematics, Bar-Ilan University, Ramat Gan 52900, Israel}
\email{machura@math.biu.ac.il}

\author{Saharon Shelah}
\address[Shelah]{Einstein Institute of Mathematics, The Hebrew University of Jerus\-alem,
Givat Ram, 91904 Jerus\-alem, Israel; and Mathematics Department,
Rutgers University, New Brunswick, NJ, USA}
\email{shelah@math.huji.ac.il}

\author{Boaz Tsaban}
\address[Tsaban]{Department of Mathematics, Bar-Ilan University, Ramat Gan 5290002, Israel; and Department of Mathematics, Weizmann Institute of Science, Rehovot 7610001, Israel}
\email{tsaban@math.biu.ac.il}
\urladdr{http://www.cs.biu.ac.il/\~{}tsaban}

\keywords{%
pseudointersection number, linear refinement number, forcing, Mathias forcing,
$\omega$-cover, $\gamma$-cover, $\tau$-cover, $\tau^*$-cover,
selection principles}

\subjclass{%
Primary: 03E17 
Secondary: 03E75 
}

\begin{document}

\begin{abstract}
The \emph{linear refinement number} $\lr$ is the minimal cardinality of a centered family in $\roth$ such that
no linearly ordered set in $(\roth,\as)$ refines this family. The \emph{linear excluded middle number} $\lx$ is
a variation of $\lr$.
We show that these numbers estimate the critical cardinalities of a number of selective covering properties.
We compare these numbers to the classic combinatorial cardinal characteristics of the continuum.
We prove that $\lr=\lx=\fd$ in all models where the continuum is at most $\aleph_2$,
and that the cofinality of $\lr$ is uncountable.
Using the method of forcing, we show that $\lr$ and $\lx$ are not provably equal to $\fd$,
and rule out several potential bounds on these numbers.
Our results solve a number of open problems.
\end{abstract}

\maketitle

\section{Overview}

\subsection{Combinatorial cardinal characteristics of the continuum}
The definitions and basic properties not included below can be found in \cite{BlassHBK}.

A family $\cF\sub\roth$ is \emph{centered} if every finite subset of $\cF$ has an infinite intersection.
For $A,B\in\roth$, $B\as A$ means that $B\sm A$ is finite.
A \emph{pseudointersection} of a family $\cF\sub\roth$ is an element $A\in\roth$ such that $A\as B$
for all $B\in\cF$.
The \emph{pseudointersection number} $\fp$ is the minimal cardinality of a centered family in $\roth$ that has no pseudointersection.

\bdfn[\cite{tautau}]
A family $\cF\sub\roth$ is \emph{linear} if it is linearly ordered by $\as$.
A family $\cG\sub\roth$ is a \emph{refinement} of a family $\cF\sub\roth$ if for each $A\in\cF$ there
is $B\in\cG$ such that $B\as A$.
The \emph{linear refinement number}
$\lr$ is the minimal cardinality of a centered family in $\roth$ that has no linear refinement.\footnote{In \cite{tautau},
the ad-hoc name $\fp^*$ is used for the linear refinement number, that we name here $\lr$.}
\edfn

A \emph{tower} is a linear subset of $\roth$ with no pseudointersection.
The \emph{tower number}  $\ft$ is the minimal cardinality of a tower.
It is immediate from the definitions that $\fp=\min\{\ft,\lr\}$.
Solving a longstanding problem, Malliaris and the second named author have recently proved that
$\fp=\ft$ \cite{MalliarisShelah}. We prove that, consistently, $\fp<\lr<\fc$.
This settles \cite[Problem 64]{tautau} (quoted in \cite[Problem 5]{ShTb768} and in \cite[Problem 11.2 (311)]{OPiT}).
Moreover, we have that
$\lr=\fd$ in all models of set theory where the continuum is at most $\aleph_2$.
One of our main results is that the cofinality of $\lr$ is uncountable.
The proof of this result uses auxiliary results of independent interest. One striking consequence is that if $\fp<\fb,$ then $\lr\le\fb$.

The number defined below is a variation of $\lr$.

\bdfn[\cite{tautau}]
For $f,g\in\NN$, let $\ev{f\le g}=\smallset{n}{f(n)\le g(n)}$.
The \emph{linear excluded middle number}
$\lx$ is the minimal cardinality of a set $\cF\sub\NN$ such that, for each $h\in\NN$,
the family $\smallset{\ev{f\le h}}{f\in\cF}$ (is either not contained in $\roth$, or) does not have a linear refinement.\footnote{In \cite{tautau}, the ad-hoc name \emph{weak excluded middle number} ($\mathfrak{wx}$)
is used for the linear excluded middle number.
Since the \emph{excluded middle number} defined in \cite{tautau} turned out equal to the classic cardinal
$\max\{\fb,\fs\}$, there is no point in preserving this name, and consequently also the name of its weaker version.
}
\edfn

If $\cF\sub\NN$ and $\card{\cF}<\lx$ then there are $h\in\NN$ and infinite subsets $A_f\as\ev{f\le h}$ such that
the family $\smallset{A_f}{f\in\cF}$ is linear, and for all $f,g\in\cF$, say such that $A_f\as A_g$, we have that $h$
excludes middles in the sense that
$$f(n)\le h(n)<g(n)$$
may hold for at most finitely many $n$ in $A_f$.

It is known that $\lr\le\lx\le\fd$ \cite{tautau} and that $\fb,\fs\le\lx$ \cite{ShTb768}.
In particular, by the above-mentioned result on $\lr$, we have that
$\lx=\fd$ whenever the continuum is at most $\aleph_2$.
In light of the results of \cite{MShT858}, Problem 57 in \cite{tautau}
asks whether $\lx=\max\{\fb,\fs\}$.
The answer, provided here, is ``No'': In the model obtained by adding $\aleph_2$ Cohen
reals to a model of \CH{}, $\fb=\fs=\aleph_1<\fd$, and thus also $\fb=\fs<\lx=\aleph_2$ in this model.
This also answers the question whether $\mathsf{wX} = \mathsf{X}$, posed in \cite{tautau} before
Problem 57, since the critical cardinalities (defined below)
of $\mathsf{wX}$ and $\mathsf{X}$ are $\lx$ and $\max\{\fb,\fs\}$,
respectively.

For $\lx$, an assertion finer than the above-mentioned one holds: If
$\fb=\lx$, then $\lx=\fd$.

We use the method of forcing (necessarily, beyond continuum of size $\aleph_2$), to show
that, consistently, $\lr,\lx<\fd$, and to rule out a number of potential upper or lower bounds on these
relatively new numbers in terms of classic combinatorial cardinal characteristics of the continuum.
We conclude by stating a number of open problems.
\subsection{Selective covering properties}

Topological properties defined by diagonalizations of open or Borel covers
have a rich history in various areas of general topology and analysis,
see \cite{LecceSurvey, KocSurv, OPiT, SakaiScheepers} for surveys on the topic
and some of its applications and open problems.

Let $X$ be an infinite topological space.
By \emph{a cover of $X$} we mean a family $\cU$ with $X\nin\cU$ and $X=\Un\cU$.
Let $\cU=\smallset{U_n}{n<\w}$ be a bijectively enumerated, countably infinite cover of $X$. We say that:
\be
\item $\cU\in\Op(X)$ if each $U_n$ is open.
\item $\cU\in\Omega(X)$ if $\cU\in\Op(X)$, and each finite subset of $X$ is contained in some $U_n$.
\item $\cU\in\Tau^*(X)$ if $\cU\in\Op(X)$, the sets
$$\set{n}{x\in U_n}\quad (\mbox{for }x\in X)$$
are infinite, and the family of these sets has a linear refinement.
\item $\cU\in\Gamma(X)$ if $\cU$ is a point-cofinite cover, that is,
each element of $X$ is a member of all but finitely many $U_n$.
\ee
We may omit the part ``$(X)$'' from these notations.

\begin{figure}[!ht]
\renewcommand{\sr}[2]{{\txt{$#1$\\$#2$}}}
{\tiny
\begin{changemargin}{-3cm}{-3cm}
\begin{center}
$\xymatrix@C=7pt@R=6pt{
&
&
& \sr{\ufin(\Op,\Ga)}{\fb}\ar[r]
& \sr{\ufin(\Op,\Tau^*)}{\lx}\ar[rr]
&
& \sr{\ufin(\Op,\Omega)}{\fd}\ar[rrrr]
&
&
&
& \sr{\sfin(\Op,\Op)}{\fd}
\\
&
&
& \sr{\sfin(\Ga,\Tau^*)}{\boxed{\lx}}\ar[rr]\ar[ur]
&
& \sr{\sfin(\Ga,\Omega)}{\fd}\ar[ur]
\\
\sr{\sone(\Ga,\Ga)}{\fb}\ar[uurrr]\ar[rr]
&
& \sr{\sone(\Ga,\Tau^*)}{\boxed{\lx}}\ar[ur]\ar[rr]
&
& \sr{\sone(\Ga,\Omega)}{\fd}\ar[ur]\ar[rr]
&
& \sr{\sone(\Ga,\Op)}{\fd}\ar[uurrrr]
\\
&
&
& \sr{\sfin(\Tau^*,\Tau^*)}{\textbf{?}}\ar'[r][rr]\ar'[u][uu]
&
& \sr{\sfin(\Tau^*,\Omega)}{\fd}\ar'[u][uu]
\\
\sr{\sone(\Tau^*,\Ga)}{\fp}\ar[rr]\ar[uu]
&
& \sr{\sone(\Tau^*,\Tau^*)}{\textbf{?}}\ar[uu]\ar[ur]\ar[rr]
&
& \sr{\sone(\Tau^*,\Omega)}{\boxed{\od}}\ar[uu]\ar[ur]\ar[rr]
&
& \sr{\sone(\Tau^*,\Op)}{\boxed{\od}}\ar[uu]
\\
&
&
& \sr{\sfin(\Omega,\Tau^*)}{\lr}\ar'[u][uu]\ar'[r][rr]
&
& \sr{\sfin(\Omega,\Omega)}{\fd}\ar'[u][uu]
\\
\sr{\sone(\Omega,\Ga)}{\fp}\ar[uu]\ar[rr]
&
& \sr{\sone(\Omega,\Tau^*)}{\boxed{\min\{\cov(\cM),\lr\}}}\ar[uu]\ar[ur]\ar[rr]
&
& \sr{\sone(\Omega,\Omega)}{\cov(\cM)}\ar[uu]\ar[ur]\ar[rr]
&
& \sr{\sone(\Op,\Op)}{\cov(\cM)}\ar[uu]
}$
\end{center}
\end{changemargin}
}
\caption{The surviving properties}\label{tau*Sch}
\end{figure}
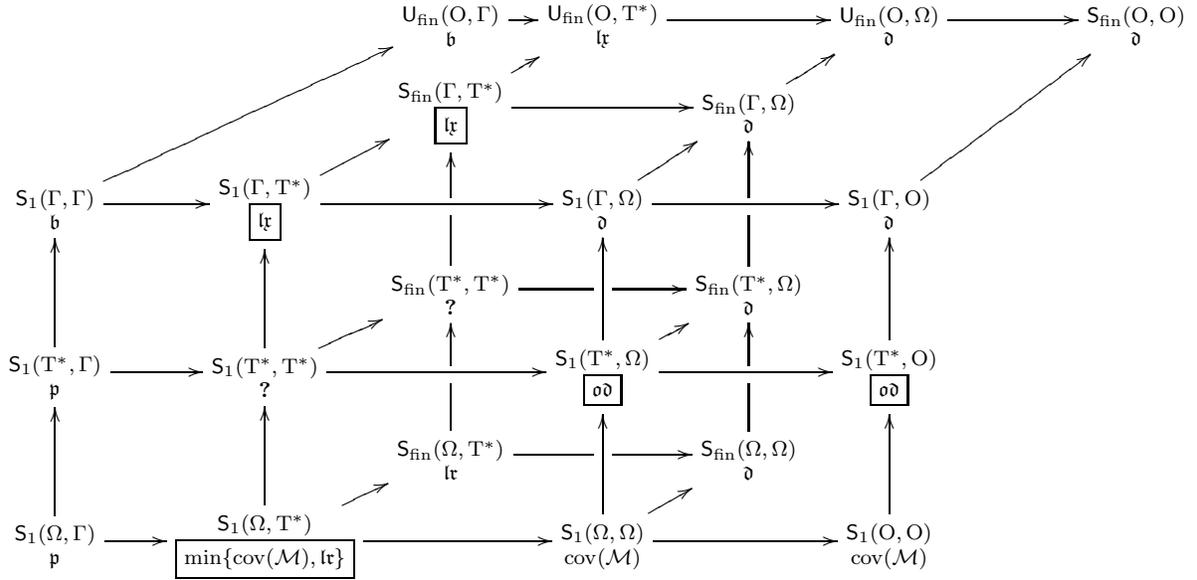

Let $\rmA$ and $\rmB$ be any of the above four types of open covers.
Scheepers \cite{coc1} introduced the following selection hypotheses that the space
$X$ may satisfy:
\bi
\item[$\sone(\rmA,\rmB)$:]
For each sequence $\smallseq{\cU_n}{n<\w}$ of members of $\rmA$,
there is a selection $\smallseq{U_n\in\cU_n}{n\in\cU_n}$ such that $\setseq{U_n}\in\rmB$.
\item[$\sfin(\rmA,\rmB)$:]
For each sequence $\smallseq{\cU_n}{n<\w}$
of members of $\rmA$, there is a selection of finite sets
$\smallseq{\cF_n\sub\cU_n}{n<\N}$ such that $\Un_{n<\N}\cF_n\in\rmB$.
\item[$\ufin(\rmA,\rmB)$:]
For each sequence $\smallseq{\cU_n}{n<\w}$ of members of $\rmA$
which do not contain a finite subcover,
there is a selection of finite sets
$\smallseq{\cF_n\sub\cU_n}{n<\N}$
such that $\setseq{\Un\cF_n}\in\rmB$.
\ei
Some of the properties are never satisfied, and many equivalences hold among
the meaningful ones. The surviving properties appear in Figure \ref{tau*Sch},
where an arrow denotes implication \cite{tautau}.
It is not known whether any implication, that does not follow from composition of existing ones,
can be added to this diagram. Several striking results concerning this problem were established by
Zdomskyy in \cite{SF2}.

Below each property $P$ in Figure \ref{tau*Sch} appears its \emph{critical cardinality}, $\non(P)$, which is the minimal
cardinality of a space $X$ not satisfying that property.\footnote{The cardinal
$\od$ was defined in \cite{MShT858}. We recall the definition in Subsection \ref{sec:ZFCSPM}, where it is needed.}
The boxed critical cardinalities, and several critical cardinalities
of properties not displayed here, are established in the present paper.

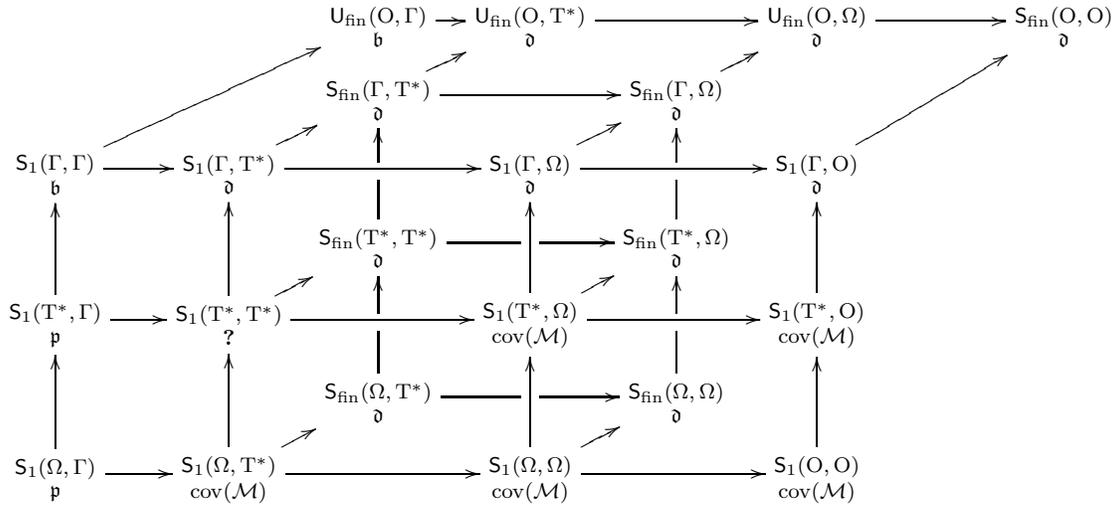
\begin{figure}[!ht]
\renewcommand{\sr}[2]{{\txt{$#1$\\$#2$}}}
{\tiny
\begin{changemargin}{-3cm}{-3cm}
\begin{center}
$\xymatrix@C=7pt@R=6pt{
&
&
& \sr{\ufin(\Op,\Ga)}{\fb}\ar[r]
& \sr{\ufin(\Op,\Tau^*)}{\fd}\ar[rr]
&
& \sr{\ufin(\Op,\Omega)}{\fd}\ar[rrrr]
&
&
&
& \sr{\sfin(\Op,\Op)}{\fd}
\\
&
&
& \sr{\sfin(\Ga,\Tau^*)}{\fd}\ar[rr]\ar[ur]
&
& \sr{\sfin(\Ga,\Omega)}{\fd}\ar[ur]
\\
\sr{\sone(\Ga,\Ga)}{\fb}\ar[uurrr]\ar[rr]
&
& \sr{\sone(\Ga,\Tau^*)}{\fd}\ar[ur]\ar[rr]
&
& \sr{\sone(\Ga,\Omega)}{\fd}\ar[ur]\ar[rr]
&
& \sr{\sone(\Ga,\Op)}{\fd}\ar[uurrrr]
\\
&
&
& \sr{\sfin(\Tau^*,\Tau^*)}{\fd}\ar'[r][rr]\ar'[u][uu]
&
& \sr{\sfin(\Tau^*,\Omega)}{\fd}\ar'[u][uu]
\\
\sr{\sone(\Tau^*,\Ga)}{\fp}\ar[rr]\ar[uu]
&
& \sr{\sone(\Tau^*,\Tau^*)}{\textbf{?}}\ar[uu]\ar[ur]\ar[rr]
&
& \sr{\sone(\Tau^*,\Omega)}{\cov(\cM)}\ar[uu]\ar[ur]\ar[rr]
&
& \sr{\sone(\Tau^*,\Op)}{\cov(\cM)}\ar[uu]
\\
&
&
& \sr{\sfin(\Omega,\Tau^*)}{\fd}\ar'[u][uu]\ar'[r][rr]
&
& \sr{\sfin(\Omega,\Omega)}{\fd}\ar'[u][uu]
\\
\sr{\sone(\Omega,\Ga)}{\fp}\ar[uu]\ar[rr]
&
& \sr{\sone(\Omega,\Tau^*)}{\cov(\cM)}\ar[uu]\ar[ur]\ar[rr]
&
& \sr{\sone(\Omega,\Omega)}{\cov(\cM)}\ar[uu]\ar[ur]\ar[rr]
&
& \sr{\sone(\Op,\Op)}{\cov(\cM)}\ar[uu]
}$
\end{center}
\end{changemargin}
}
\caption{The critical cardinalities in models of $\fc\le\aleph_2$}\label{fig:daleph2}
\end{figure}

Putting the mentioned results together, we have that in models where the continuum (or just $\fd$) is at most $\aleph_2$,
all but one of the critical cardinalities of the studied properties are determined in terms of classic combinatorial cardinal characteristics of the continuum, see Figure \ref{fig:daleph2}.

These results fix, in particular, an erroneous assertion made in \cite[Theorem 7.20]{tautau} without proof,
namely, that the critical cardinality of $\sone(\Om,\Tau^*)$ is $\lr$. As shown in the diagram, the correct
critical cardinality is $\min\{\cov(\cM),\lr\}$. By the above-mentioned results,
the inequality $\cov(\cM)<\lr$ holds in all models of $\cov(\cM)<\fd=\aleph_2$;
in particular in the standard Laver, Mathias, and Miller models (see \cite{BlassHBK}).

\section{Results in ZFC}

\subsection{Combinatorial cardinal characteristics of the continuum}
All filters in this paper are on $\N$, and \emph{are assumed to contain all cofinite subsets} of $\N$.
The \emph{character} of a filter $\cF$ is the minimal cardinality of a base for $\cF$, that is, a set $\cB\sub\cF$ such that each element of $\cF$ contains some element of $\cB$,
or equivalently, the minimal cardinality of a subset $\cB$ of $\cF$ generating $\cF$ as a filter.
Let $\cF$ be a filter. A set $P\sub\N$ is \emph{$\cF$-positive} if $P\cap A$ is infinite for all $A\in\cF$, in other words, $\cF$ can be extended to a filter containing $P$.

\blem\label{lem:filext}
Let $\kappa$ be an infinite cardinal such that, for each filter $\cF$ of character $\le\kappa$, every linear
subset of $\cF$ of cardinality $<\kappa$ has an $\cF$-positive pseudointersection.
Then $\kappa<\lr$.
\elem
\bpf
Let $\smallset{A_{\alpha}}{\alpha < \kappa}$ be centered, and $\cF$ be the filter generated by $\smallset{A_{\alpha}}{\alpha < \kappa}$.
We construct a linear refinement $\smallset{A_\alpha^-}{\alpha<\kappa}$ of $\smallset{A_{\alpha}}{\alpha < \kappa}$ by induction on $\alpha$.

Let $A_0^-=A_0$. For $\alpha>0$ we assume, inductively, that
$\smallset{A_{\beta}^-}{\beta < \alpha}$ is linear and that $\cF\cup\smallset{A_{\beta}^-}{\beta < \alpha}$ is centered.
Let $\cF_\alpha$ be the filter generated by $\cF\cup\smallset{A_{\beta}^-}{\beta < \alpha}$.
Let $P$ be an $\cF_\alpha$-positive pseudointersection of $\smallset{A_{\beta}^-}{\beta < \alpha}$.
Take $A_\alpha^-=P\cap A_\alpha$.
As $\cF$ is a filter, $A_\alpha^-$ is $\cA$-positive.
\epf

In the following proof, we use that $\lr\le\fd$ \cite{tautau}.
Theorem \ref{thm:selOmT} improves upon this inequality.

\bthm\label{thm:aleph1}
If $\lr=\aleph_1$, then $\fd=\aleph_1$.
\ethm
\bpf
Assume that $\fd>\aleph_1$. We will prove, using Lemma \ref{lem:filext}, that $\lr>\aleph_1$.
Let $\cF$ be a filter of character $\le\aleph_1$, and fix a base $\smallset{B_\alpha}{\alpha<\aleph_1}$ of $\cF$.
Let $\smallset{A_n}{n<\N}$ be a linear subset of $\cF$.
We prove that $\smallset{A_n}{n<\N}$ has an $\cF$-positive pseudointersection.
We may assume that $A_{n+1}\sub A_n$ for all $n$.

Let $\alpha<\aleph_1$. For each $n$, as $B_\alpha\cap A_n\in\cF$,
we can pick an element
$$f_\alpha(n)\in B_\alpha\cap A_n$$
such that the function $f_\alpha$ is strictly increasing.
As $\fd>\aleph_1$, there is $g\in\NN$ such that, for each $\alpha<\aleph_1$, $f_\alpha(n)<g(n)$ for infinitely many $n$.
Let
$$P=\Un_{n<\N}A_n\cap \bintvl{0}{g(n)}.$$
For each $n$, $P\sm A_n\sub \Un_{k<n} \bintvl{0}{g(k)}$, and thus $P\as A_n$.
For each $\alpha<\aleph_1$ and each $n$ with $f_\alpha(n)<g(n)$,
$$f_\alpha(n)\in B_\alpha\cap A_n\cap \bintvl{0}{g(n)}\sub P.$$
As $f_\alpha$ is strictly increasing, $B_\alpha\cap P$ is infinite.
Thus, $P$ is $\cF$-positive.
\epf

As $\lr\le\fd$ \cite{tautau}, it follows from Theorem \ref{thm:aleph1} that, if $\fd\le\aleph_2$, then $\lr=\fd$.
Thus, a large family of results about combinatorial cardinal characteristics of the continuum in models of $\fc=\aleph_2$ (see Table 4 in \cite{BlassHBK})
are applicable. For example, we have the following consequences.

\bcor\label{cor:conslr}
\mbox{}
\be
\itm For each cardinal $\mathfrak x$ among $\mathfrak r$, $\mathfrak u$, $\mathfrak a$, $\cov(\cN)$, $\non(\cN)$, and $\non(\cM)$,
it is consistent that $\mathfrak x<\lr$, and it is consistent that $\lr<\mathfrak x$.

\itm For each cardinal $\mathfrak x$ among $\mathfrak p$, $\mathfrak h$, $\mathfrak s$, $\mathfrak g$, $\mathfrak e$, $\mathfrak b$, $\add(\cN)$, $\add(\cM)$, and $\cov(\cM)$, it is consistent that $\mathfrak x<\lr$.

\itm For each cardinal $\mathfrak x$ among $\mathfrak i$, $\cof(\cM)$, and $\cof(\cN)$,
it is consistent that $\lr<\mathfrak x$.
\qed
\ee
\ecor

In Subsection \ref{sec:lrLTcovM} we show that, consistently, $\lr<\cov(\cM)$. In particular, $\lr<\fd$ is consistent.

A \emph{tower of height $\kappa$} is a set $\smallset{T_\alpha}{\alpha<\kappa}\sub\roth$ that is $\as$-decreasing
with $\alpha$ and has no pseudointersection. There is no tower of height smaller than $\fp$, and by the Malliaris--Shelah Theorem, $\fp$ is the minimal height of a tower.

\blem\label{tref}
Let $\cF\sub \roth$ be a centered family of cardinality smaller than $\lr$.
Then either $\cF$ has a pseudointersection, or $\cF$ is refined by a tower of height $\fp$.
\elem
\bpf
If $\fp=\lr$, then $\cF$ has a pseudointersection, and we are done.

Assume that $\fp<\lr$.
Let $\smallset{P_\alpha}{\alpha<\fp}\sub \roth$ be a centered family with no pseudointersection.
Set
$$\cB=\set{A\x P_\alpha}{A\in\cF,\alpha<\fp}\cup \set{\smallset{(n,m)}{k\le n,m}}{k\in \w}.$$
Then $\cB$ is a centered family of cardinality
less than $\lr$. Let $\cR=\smallset{R_\alpha}{\alpha<\kappa}\sub[\w\x\w]^\w$
be a $\as$-decreasing linear refinement of $\cB$, with $\kappa$ regular.

Let $\pi_0$ and $\pi_1$ be the projections of $\w\x\w$ on the first and
second coordinates, respectively. For each pseudointersection $R$ of
the family $\smallset{\smallset{(n,m)}{k\le n,m}}{k\in \w}$, the sets
$\pi_0(R)$ and $\pi_1(R)$ are both infinite.
Moreover, if $R\as A\x B$ then $\pi_0(R)\as A$ and $\pi_1(R)\as B$.

If $\kappa<\fp$, then $\cR$ has a pseudointersection $R$.
By the above paragraph, the set $A:=\pi_0(R)$ is infinite, and
is a pseudointersection of $\cF$, as required.

Next, assume that $\fp\le\kappa$. For each $k$, fix $\alpha_k$ such that
$R_{\alpha_k}\as \smallset{(n,m)}{k\le n,m}$. As $\kappa$ is uncountable and regular, we have that $\alpha :=\sup_k\alpha_k<\kappa$. Removing the first $\alpha$ members of $\cR$, we may assume that every member of $\cR$ is a pseudointersection of the family $\smallset{\smallset{(n,m)}{k\le n,m}}{k\in \w}$, and, consequently, that the sets $\pi_0(R)$ and $\pi_1(R)$ are infinite for each $R\in\cR$.
It follows that the families $\smallset{\pi_0(R_\alpha)}{\alpha<\kappa}$
and $\smallset{\pi_1(R_\alpha)}{\alpha<\kappa}$ are linear refinements of
the families $\cF$ and $\smallset{P_\alpha}{\alpha<\fp}$, respectively.
In particular, if $\kappa=\fp$, then we are done.

It remains to prove that the case $\kappa>\fp$ is impossible.
Assume otherwise. For each $\alpha<\fp$, fix $\beta_\alpha<\kappa$
such that $\pi_1(R_{\beta_\alpha})\as P_\alpha$.
As $\kappa$ is regular, we have that $\beta:=\sup_{\alpha<\ft}\beta_\alpha<\kappa$,
and $\pi_1(R_\beta)$ is a pseudointersection of the family $\smallset{P_\alpha}{\alpha<\fp}$; a contradiction.
\epf

\blem[Folklore]
If $\fb<\fd$ then there is a tower of height $\fb$.
\elem
\bpf
Let $\smallset{f_\alpha}{\alpha<\fb}\sub\NN$ be a $\fb$-scale, that is, an unbounded set
where each $f_\alpha$ is an increasing member of $\NN$ and the sequence $f_\alpha$ is $\le^*$-increasing with $\alpha$.
Let $h\in\NN$ witness that this family is not dominating.
Then $\smallset{\ev{f_\alpha\le h}}{\alpha<\fb}$ is a tower, for if $P$ is a pseudointersection, then $\smallset{f_\alpha\rest P}{\alpha<\fb}$ is bounded by $h\rest P$.
\epf

\bthm\label{surprise}
If $\fp<\fb$ then $\lr\le\fb$.
\ethm
\bpf
Assume that $\fb<\lr$. Then, as $\lr\le\fd$, we have that $\fb<\fd$ and there is a tower $\smallset{T_\alpha}{\alpha<\fb}$
of height $\fb$.
By Lemma \ref{tref}, this tower is refined by a tower  $\smallset{P_\alpha}{\alpha<\fp}$.
Assume that $\fp<\fb$.
For each $\alpha<\fp$, fix $\beta_\alpha<\fb$ with
$T_{\beta_\alpha}\not\as P_\alpha$. As $\fb$ is regular, $\beta:=\sup_{\alpha<\ft}\beta_\alpha<\fb$.
Then $T_\beta$ is not refined by any $P_\alpha$; a contradiction.
\epf

The argument in the last proof shows the following.

\bcor
Each tower of regular height smaller than $\lr$ must be of height $\fp$.\qed
\ecor

A family $\cF\sub\NN$ is \emph{$\kappa$-bounded} if there is a family $\cG\sub\NN$ of cardinality $\kappa$ such that each member of $\cF$ is dominated by some member of $\cG$.

\blem\label{tbdd}
Let $\cF\sub\NN$.
\be
\item If $\card{\cF}<\lr$, then $\cF$ is  $\fp$-bounded.
\item If $\cof(\lr)\le \fp$ and $\card{\cF}=\lr$, then $\cF$ is  $\fp$-bounded.
\ee
\elem
\bpf
(1) Let $\cF\sub\NN$.
We may assume that each member of $\cF$ is an increasing function.

Assume that $\card{\cF}<\lr$.
For each $f\in\cF$, let
$$A_f=\set{(n,m)}{f(n)\le m}\sub\w\x\w.$$
The family
$$\set{A_f}{f\in\cF}\cup\set{\set{(n,m)}{n>k}}{k\in\N}$$
is centered.

Assume that this family has a pseudointersection $A$. As $A$ is a
pseudointersection of $\smallset{\smallset{(n,m)}{n>k}}{k\in\N}$,
infinitely many columns $A\cap\{n\}\x\w$ of $A$ (for $n<\N$) are nonempty,
and all columns of $A$ are finite.
For each $n$, define $g_A(n)$ as follows:
Let $n'\ge n$ be minimal with the column $A\cap(\{n'\}\x\w)$ nonempty, and let $g_A(n)$ be minimal such that $(n',g_A(n))$ is in that column.
For each $f\in\cF$, as $A\as A_f$ and $f$ is increasing, we have that $f\le^* g_A$. Thus, $\cF$ is bounded, and we are done.

Next, assume that our family does not have a pseudointersection.
By Lemma \ref{tref}, some tower
$\smallset{R_\alpha}{\alpha<\fp}$ linearly refines our family.
As $\fp$ is regular, by removing an initial segment of indices we may assume that each $R_\alpha$ is a pseudointersection of
$\smallset{\smallset{(n,m)}{n>k}}{k\in\N}$.
Thus, we can define functions $g_{R_\alpha}$ for $\alpha<\fp$ as in the previous paragraph. As above, for each $f\in \cF$, if $\alpha<\fp$ is such that $R_\alpha\as A_f$, then $f\le^* g_{R_\alpha}$. This shows that $\cF$ is $\fp$-bounded.

(2) Assume that $\card{\cF}=\lr$.
Represent $\cF=\Un_{\alpha<\cof(\lr)}\cF_\alpha$, with $\card{\cF_\alpha}<\lr$ for each $\alpha$. Then every $\cF_\alpha$ is $\fp$-bounded.
As $\cof(\lr)\cdot\fp=\fp$, $\cF$ is $\fp$-bounded.
\epf

\bthm\label{thm:cflr}
The cofinality of $\lr$ is uncountable.
\ethm
\bpf
As $\fp$ is regular, we have that $\lr$ is regular if $\lr=\fp$.

Assume that $\fp<\lr$. Towards a contradiction, assume that $\cof(\lr)=\alephes$.
Let
$$\cF=\set{A_\alpha}{\alpha<\lr}\sub\roth$$
be a centered family. We will prove that $\cF$ has a linear refinement. Represent $\cF=\Un_n\cF_n$ with $\cF_n\sub\cF_{n+1}$ and $\card{\cF_n}<\lr$ for all $n$.
By thinning out the sequence $\smallseq{\cF_n}{n<\N}$, we may assume that
each $\cF_n$ has a pseudointersection, or no $\cF_n$ has a pseudointersection.

Consider first the former case. For each $n$, let $R_n$ be a pseudointersection
of $\cF_n$. For each $A\in\cF$, let $k$ be the first with $A\in\cF_k$. For $n<k$ let $f_A(n)=0$, and for $n\ge k$ let
$$f_A(n)=\min\set{m}{R_n\sm m\sub A}.$$
By Lemma \ref{tbdd}, the family $\smallset{f_A}{A\in\cF}$ is $\fp$-bounded. Let $\cG\sub\NN$ be a witness for that. For each $g\in\cG$ and each $k$, let
$$U_{g,k}=\Un_{n\ge k}R_n\sm g(n).$$
The family $\smallset{U_{g,k}}{g\in\cG, k\in\N}$ is centered.
Indeed, for $k_1,\dots,k_l\in\N$ and $g_1,\dots,g_l\in\NN$, let $n=\max\{k_1,\dots,k_l\}$ and $m=\max\{g_1(n),\dots,g_l(n)\}$.
Then $n\ge k_1,\dots,k_l$ and $R_n\sm m\sub U_{g_1,k_1}\cap \cdots \cap U_{g_l,k_l}$.
Since the cardinality of this family is at most $\fp<\lr$,
it has a linear refinement $\cR$.
Let $A\in\cF$, and let $g\in\cG$ be such that $f_A\le^* g$. Fix $k$ such that $f_A(n)\le g(n)$ for all $n\ge k$. Then $U_{g,k}\sub A$. Thus, $\cR$ is also a linear refinement of $\cF$.

It remains to consider the case where no $\cF_n$ has a pseudointersection.
This is done by slightly extending the previous argument.
By Lemma \ref{tref},
for each $n$, there is a tower $\smallset{T^n_\alpha}{\alpha<\fp}$ that linearly refines $\cF_n$.
Fix $A\in\cF$, and let $k$ be the first with $A\in\cF_k$. For $n<k$ let
$\alpha_n=0$, and for $n\ge k$ let $\alpha_n<\fp$ be the first with
$T^n_{\alpha_n}\as A$. As $\fp$ is regular, the ordinal $\alpha(A):=\sup_n\alpha_n$ is smaller than $\fp$. Then
$$T^n_{\alpha(A)}\as A$$
for all but finitely many $n$.
For $n<k$ let $f_A(n)=0$, and for $n\ge k$ let
$$f_A(n)=\min\set{m}{T^n_{\alpha(A)}\sm m\sub A}.$$

By Lemma \ref{tbdd}, the family $\smallset{f_A}{A\in\cF}$ is $\fp$-bounded. Let $\cG\sub\NN$ be a witness for that. For each $g\in\cG$, $\alpha<\fp$ and $k\in\N$, let
$$U_{g,\alpha,k}=\Un_{n\ge k}T^n_{\alpha}\sm g(n).$$
The family $\smallset{U_{g,\alpha,k}}{g\in\cG, \alpha<\fp, k\in\N}$ is centered, and has cardinality $\fp<\lr$. Thus, it has a linear refinement $\cR$.
Let $A\in\cF$, and let $g\in\cG$ be such that $f_A\le^* g$. Fix $k$ such that $f_A(n)\le g(n)$ for all $n\ge k$. Then $U_{g,\alpha(A),k}\sub A$. Thus, $\cR$ is also a linear refinement of $\cF$.
\epf

We conclude this subsection with a result on $\lx$ that is analogous to Theorem \ref{thm:aleph1}.
Recall from Figure \ref{tau*Sch} that $\fb\le\lx\le\fd$.

\bthm
If $\lx=\fb$ then $\fd=\fb$.
\ethm
\bpf
Assume that $\fb<\fd$.
Let $\smallset{f_\alpha}{\alpha<\fb}\sub\NN$.
We will find a function $h\in\NN$ and a linear refinement of the family
$\smallset{\ev{f_\alpha\le h}}{\alpha<\fb}$.

For each $\alpha<\fb$, let $g_\alpha$ be a $\le^*$-bound of
$\{f_\alpha\}\cup\smallset{g_\beta}{\beta<\alpha}$.
Let $h\in\NN$ witness that $\smallset{g_\alpha}{\alpha<\fb}$ is not
dominating. Then
$\smallset{\ev{g_\alpha\le h}}{\alpha<\fb}$ is a linear refinement of $\smallset{f_\alpha}{\alpha<\fb}$.
\epf

As $\lr\le\lx\le\fd$, Corollary \ref{cor:conslr} holds for $\lx$ as well.
In Section \ref{sec:lrLTcovM} we show that, consistently, $\lx<\fd$.

\subsection{Selective covering properties}\label{sec:ZFCSPM}

For a topological space $X$, let $\Tau(X)$ denote the family of all
open covers $\smallset{U_n}{n<\w}$ of $X$ such that the sets
$\smallset{n}{x\in U_n}$ (for $x\in X$) are infinite, and the family of these sets is linear.
The first result of this section solves one of the first problems concerning this type of covers \cite[Problem 10]{tautau}
(quoted in \cite[Problem 7.2]{OPiT}).

Let $\smallbinom{\rmA}{\rmB}$ denote the property that every element of $\rmA$ contains an element of $\rmB$.

\bthm\label{thm:selOmT}
Let $\rmA\sub\Tau^*$. Then  $\smallbinom{\Om}{\rmA}=\sfin(\Om,\rmA)$. In particular:
\be
\itm $\smallbinom{\Om}{\Tau}=\sfin(\Om,\Tau)$; and
\itm $\smallbinom{\Om}{\Tau^*}=\sfin(\Om,\Tau^*)$.
\ee
\ethm
\bpf
It suffices to prove that $\smallbinom{\Om}{\Tau^*}$ implies $\sfin(\Om,\Om)$.

Assume that $\smallset{U^n_m}{m\in\w}\in\Om(X)$ for each $n<\N$.
Fix distinct elements $x_n\in X$ for $n<\N$.
Then
$$\cU:=\set{U^n_m\sm\{x_n\}}{n,m\in\N}\in\Om(X).$$
Let $\cV\sub\cU$ be such that $\cV\in\Tau^*(X)$.
Enumerate $\cV=\smallset{V_n}{n<\N}$.
For $x\in X$, let $\cV(x)=\smallset{n}{x\in V_n}$.
By the definition of $\Tau^*$, the family $\smallset{\cV(x)}{x\in X}$ has a linear refinement $\cR$.

There is an pseudointersection $P$ of the family $\smallset{\cV(x_n)}{n<\N}$
such that, for each finite $\cF\sub\cR$, $P\cap\bigcap\cF$ is infinite.
Indeed, if $\cR$ has a pseudointersection then we can take $P$ to be this pseudointersection.
And if not, then by thinning $\cR$ out, we may assume that
$\cR=\smallset{R_\alpha}{\alpha<\kappa}$ is a tower of regular uncountable height $\kappa$.
For each $n$, let $\alpha_n<\kappa$ be with $R_{\alpha_n}\as\cV(x_n)$. Let $\alpha=\sup_n\alpha_n$, and
take $P=R_\alpha$.

Let $\cW=\smallset{V_k}{k\in P}$.
Fix $n$. As $P\as\cV(x_n)$, we have that $x_n\in V_k$ for all but finitely many $k\in P$.
Thus, the set $\cW\cap \smallset{U^n_m\sm\{x_n\}}{m\in\N}$ is finite.
Let $F_n\sub\N$ be a finite (possibly empty) set such that
$$\cW\cap \set{U^n_m\sm\{x_n\}}{m\in\N}=\set{U^n_m\sm\{x_n\}}{m\in F_n}.$$
Let $Y$ be a finite subset of $X$. Then the set $P\cap \bigcap_{y\in Y}\cV(y)$ is infinite, and
for each $k$ in this set, $Y\sub V_k$.
Thus, $\cW\in\Om(X)$.
As $\cW\in\Om(X)$, the family $\Un_n\smallset{U^n_m}{m\in F_n}$ is in $\Om(X)$, too.
\epf

\bthm\label{nonGaTau*}
The critical cardinalities of $\sone(\Ga,\Tau^*)$ and $\sfin(\Ga,\Tau^*)$ are both $\lx$.
\ethm
\bpf
As the critical cardinality of $\ufin(\Op,\Tau^*)$ is $\lx$ \cite{tautau} and the implications
$$\sone(\Ga,\Tau^*)\longrightarrow\sfin(\Ga,\Tau^*)\longrightarrow\ufin(\Op,\Tau^*)$$
hold, it suffices to prove that every topological space of cardinality smaller than $\lx$ satisfies $\sone(\Ga,\Tau^*)$.

Let $X$ be a topological space with $\card{X}<\lx$.
Assume that, for each $n$, $\smallset{U^n_m}{m<\w}$ is a point-cofinite cover of $X$.
For each $x\in X$, define $f_x\in\NN$ by
$$f_x(n)=\min\set{m}{\forall k\ge m,\ x\in U^n_m}.$$
As $\card{X}<\lx$, there are $h\in\NN$ and infinite subsets
$$A_x\sub \ev{f_x\le h}\quad (x\in X)$$
such that $\smallset{A_x}{x\in X}$ is linear.
Then
$\smallset{U^n_{h(n)}}{n<\w}\in\Tau^*(X)$.
Indeed, for each $x\in X$,
$$A_x\sub \ev{f_x\le h}\sub \set{n}{x\in U^n_{h(n)}},$$
and the family $\smallset{A_x}{x\in X}$ is linear.
\epf

\bthm\label{nonOmTau*}
The critical cardinality of $\sone(\Om,\Tau^*)$ is $\min\{\cov(\cM),\lr\}$.
\ethm
\bpf
Notice that
$$\sone(\Om,\Tau^*)=\sone(\Om,\Om)\cap \binom{\Om}{\Tau^*}.$$
It follows that
$$\non\bigl(\sone(\Om,\Tau^*)\bigr)=
\min\bigl\{
\non\bigl(\sone(\Om,\Om)\bigr),
\non\bigl(\smallbinom{\Om}{\Tau^*}\bigr)
\bigr\}.$$
By the definitions of $\Om$ and $\Tau^*$, the critical cardinality of $\smallbinom{\Om}{\Tau^*}$ is $\lr$ \cite{tautau}.
It is known that $\non(\sone(\Om,\Om))=\cov(\cM)$.
\epf

\bthm
$\min\{\cov(\cM),\lr\},\min\{\fb,\fs\}\le\non(\sfin(\Tau^*,\Tau^*))\le\lx$.
\ethm
\bpf
By Theorem \ref{nonGaTau*}, as $\sfin(\Tau^*,\Tau^*)$ implies $\sfin(\Ga,\Tau^*)$,
we have that
$\non(\sfin(\Tau^*,\Tau^*))\le\lx$.
By Theorem \ref{nonOmTau*}, as $\sfin(\Om,\Tau^*)$ implies $\sfin(\Tau^*,\Tau^*)$,
we have that
$\min\{\cov(\cM),\lr\}\le\non(\sfin(\Tau^*,\Tau^*))$.
It remains to prove that $\min\{\fb,\fs\}\le\non(\sfin(\Tau^*,\Tau^*))$.
This is proved as in the proof of Lemma 3.4 of \cite{MShT858}.
For the reader's convenience, we provide a complete argument.

Let $X$ be a topological space with $\card{X}<\min\{\fb,\fs\}$.
Assume that, for each $n$, $\smallset{U^n_m}{m<\w}\in\Tau^*(X)$.
For each $n$, let
$$A_x(n)\sub\set{m}{x\in U^n_m}\quad (x\in X)$$
be a linear family.
For $x,y\in X$, let
$$B_{x,y}=\set{n}{A_x(n)\as A_y(n)}.$$
As $\card{X}<\fs$, there is $S\in\roth$ that is not split by any $B_{x,y}$.
As $B_{x,y}\cup B_{y,x}=\w$, $S\as B_{x,y}$ or $S\as B_{y,x}$ for all $x,y$.
For $x,y\in X$ define $g_{x,y}\in\NN$
by:
$$g_{x,y}(n) = \begin{cases}
\min\set{k}{A_x(n)\sm k\sub A_y(n)\sm k} & n \in B_{x,y}\sm B_{y,x}\\
\min\set{k}{A_y(n)\sm k\sub A_x(n)\sm k} & n \in B_{y,x}\sm B_{x,y}\\
\min\set{k}{A_x(n)\sm k=A_y(n)\sm k} & n \in B_{x,y}\cap B_{y,x}
\end{cases}$$
Since $\card{X}<\fb$, there exists $g_0\in\NN$
which dominates all of the functions $g_{x,y}$, $x,y\in X$.
For each $x\in X$, define $g_x\in\NN$ by
$$g_x(n)=\min A_x(n)\sm g_0(n).$$
Choose $g_1\in\NN$ which dominates the functions
$g_x$ (for $x\in X$). Here too, this is possible since $\card{X}<\fb$.
For each $n\in S$, let
$$\cF_n = \{U^n_{g_0(n)},\dots, U^n_{g_1(n)}\}.$$
For $n\nin S$ let $\cF_n=\emptyset$.
Let
$$\cU=\Un_{n\in S}\cF_n=\set{U^n_m}{n\in S,\ g_0(n)\le m\le g_1(n)}.$$
We claim that $\cU\in\Tau^*(X)$.
For each $x\in X$ let
$$\cU_x=\set{U^n_m}{n\in S,\ g_0(n)\le m\le g_1(n),\ m\in A_x(n)}\sub\set{U\in\cU}{x\in\cU}.$$
We may assume that the sets $U^n_m$ are distinct for distinct pairs $(n,m)$.
For all but finitely many $n\in S$,
$m:=g_x(n)\in A_x(n)$ and $g_0(n)\le g_x(n)\le g_1(n)$, so $x\in U^n_{m}\in \cU_x$.
Thus, $\cU_x$ is an infinite subset of $\cU$.
It remains to show that the family $\smallset{\cU_x}{x\in X}$ is linear.

Let $x,y\in X$. Without loss of generality, $S\as B_{x,y}$.
We will show that $\cU_x\as\cU_y$.
For all but finitely many $n\in S$:
$g_{x,y}(n)\le g_0(n)$. For each $U^n_m\in \cU_x$, $g_0(n)\le m\in A_x(n)$, and thus
$g_{x,y}(n)\le m$. As $n\in B_{x,y}$ and $m\in A_x(n)$, we have that $m\in A_y(n)$.
Thus, $U^n_m\in\cU_y$.
\epf

\bdfn[\cite{MShT858}]\label{dfn:od}
$\od$ is the minimal cardinality of a
family $\cA\sub (\roth)^\N$ such that:
\be
\item For each $n$, $\smallset{A(n)}{A\in\cA}$ is linear.
\item There is no $g\in\NN$ such that, for each $A\in\cA$,
$g(n)\in A(n)$ for some $n$.
\ee
\edfn
$\cov(\cM)\le\od$, and equality holds if the continuum is at most $\aleph_2$
\cite{MShT858}.

\bthm
The critical cardinalities of $\sone(\Tau^*,\Om)$ and of $\sone(\Tau^*,\Op)$ are both $\od$.
\ethm
\bpf
As $\Tau\sub\Tau^*$,
$$\sone(\Tau^*,\Om)\longrightarrow\sone(\Tau^*,\Op)\longrightarrow \sone(\Tau,\Op).$$
In \cite{MShT858} it is proved that $\non(\sone(\Tau,\Op))=\od$.
It remains to prove that $\od\le\sone(\Tau^*,\Om)$.

Let $X$ be a topological space with $\card{X}<\od$.
Assume that, for each $n$, $\smallset{U^n_m}{m<\w}\in\Tau^*(X)$.
Fix $n$.
By the definition of $\Tau^*$, there are sets
$$A_x(n)\sub\set{m}{x\in U^n_m}$$
such that $\smallset{A_x(n)}{x\in X}\sub\roth$ and is linear.
For each finite $F\sub X$, let
$$A_F(n)=\bigcap_{x\in F}A_x(n).$$
Then the family
$$\set{A_F(n)}{F\in [X]^{<\w}}\sub\roth$$
is linear.
As $\card{X}<\od$, there is $g\in\NN$ such that, for each finite $F\sub X$, there
is $n$ with $g(n)\in A_F(n)$.
Then $\smallset{U^n_{g(n)}}{n<\w}\in\Om(X)$.
\epf

Recall that $\fp=\ft$ \cite{MalliarisShelah}.

\bthm\label{t*t}
The critical cardinality of $\smallbinom{\Tau^*}{\Tau}$ is $\ft$.
\ethm
\bpf
We use the method of the proof of \cite[Theorem 3]{ShTb768}.

$(\ge)$ $\smallbinom{\Tau^*}{\Ga}$, which implies $\smallbinom{\Tau^*}{\Tau}$, has critical cardinality $\ft$.

$(\le)$ Consider $P(\w)$ with the Cantor space topology and the open sets
$$U_n=\set{A\in P(\w)}{n\in A}.$$
For a family $\cA\sub \roth$, viewed as a subspace of $P(\w)$:
\be
\item $\smallset{U_n}{n<\w}\in\Tau^*(\cA)$ if and only if $\cA$ has a linear refinement.
\item $\smallset{U_n}{n<\w}\in\Tau(\cA)$ if and only if $\cA$ is linear.
\item $\smallset{U_n}{n<\w}$ contains an element of $\Tau(\cA)$ if and only if
there is $I\in\roth$ such that $\smallset{A\cap I}{A\in\cA}$ is a linear subset of $\roth$.
\ee
We construct a family $\cA\sub\roth$ of cardinality $\ft$, such that
$\cA$ has a linear refinement, but for each $I\in\roth$, the family $\smallset{A\cap I}{A\in\cA}$ is nonlinear.

Let $\cF\sub\roth$ be a tower of cardinality $\ft$.
Let $\cB$ be the boolean subalgebra of $P(\N)$ generated by $\cF$.
Then $\card{\cB}=\ft$.
Let
$$\cA=\set{B\in\cB}{\exists A\in\cF,\ A\sub B}.$$
Then $\cF$ is a linear refinement of $\cA$.

Towards a contradiction, assume that there is $I\in\roth$ such that
$\smallset{A\cap I}{A\in\cA}$ is a linear subset of $\roth$.
As $\smallset{A\cap I}{A\in\cA}$ refines $\cA$, it has no pseudointersection.
Fix an element $D_0\in\cA$.
There exist:
\be
\item An element $D_1\in\cA$
such that $D_1\cap I\subset^* D_0\cap I$ (i.e., such that $D_0\cap I\sm D_1$ is infinite);
and
\item An element $D_2\in\cA$ such that $D_2\cap I\subset^* D_1\cap I$.
\ee
Then the sets $(D_2\cup(D_0\sm D_1))\cap I$ and
$D_1\cap I$, both elements of $\cA$,
contain the infinite sets
$(D_0\cap I)\sm (D_1\cap I)$ and $(D_1\cap I)\sm (D_2\cap I)$,
respectively, and thus are not $\as$-comparable, a contradiction.
\epf

\bcor
The critical cardinalities of $\sone(\Tau^*,\Tau)$ and $\sfin(\Tau^*,\Tau)$ are both $\ft$.
\ecor
\bpf
As
$$\sone(\Tau^*,\Tau)=\binom{\Tau^*}{\Tau}\cap\sone(\Tau,\Tau)$$
and $\non(\sone(\Tau,\Tau))=\ft$ \cite{MShT858},
we have by Theorem \ref{t*t} that $\non(\sone(\Tau^*,\Tau))=\ft$.
Thus, by the implications
$$\sone(\Tau^*,\Tau)\longrightarrow\sfin(\Tau^*,\Tau)\longrightarrow\binom{\Tau^*}{\Tau}$$
and Theorem \ref{t*t}, $\non(\sfin(\Tau^*,\Tau))=\ft$.
\epf

\section{Consistency results}

\subsection{A model for $\lx < \fd$}\label{sec:lrLTcovM}

For a cardinal $\lambda$, let $\bbC_\lambda$ be the forcing notion adding $\lambda$ Cohen reals.

\bthm\label{thm:Cohen}
Let $\mu=\fc$ and $\lambda > \mu^+$. Then
$$\forces_{\bbC_\lambda} \lx \leq \mu^+<\lambda\le\cov(\cM).$$
\ethm
\bpf
Let $\bbC_\lambda= \Fn( \lambda \times \omega, \omega )$ and let $c_{\alpha}$ be the $\alpha$-th Cohen real added by $\bbC_\lambda$.
For $p\in \bbC_\lambda$, let
$\supp (p)=\smallset{\beta}{\dom (p)\cap (\{\beta \} \times \omega) \ne \emptyset}$.
For $\beta \in \supp(p)$, let $p(\beta) $ be the partial function from $\omega$ to $\omega$ defined by
$p(\beta)(n) = p(\beta,n)$.
Thus, if $(\beta, n) \in \dom(p)$ and $p(\beta)(n)=m$, then $p \forces \name{c}_{\beta}(n)= m$.

We claim that the set $\smallset{c_{\alpha}}{\alpha < \mu^{+}}$ witnesses that $\lx \leq \mu^{+}$.
Towards a contradiction, assume that there are:
A condition $p\in \bbC_\lambda$, a name $\name{h}$ for a function in $\NN$,
and names $\name{A}_{\alpha}$ (for $\alpha < \mu^{+}$) of infinite subsets of $\omega$ such that
\be
\item $p\forces \name{A}_{\alpha} \as \smallset{n \in \omega}{\name{c}_{\alpha}(n) \leq \name{h}(n)} \mbox{ and } \name{A}_{\alpha} \mbox{ is infinite}$;
\item For all $\alpha$ and $\beta$, $p\forces \name{A}_{\alpha} \sub \name{A}_{\beta} \mbox{ or } \name{A}_{\beta} \sub \name{A}_{\alpha}$.
\ee
Fix $U_h \in [\lambda]^\alephes$ and a Borel function
$b_h\colon (\NN)^{U_h} \to \NN$, coded in the ground model,
such that
$$ p\forces \name{h}= b_h ( \seq{\name{c}_{\beta}}{\beta \in U_h}).$$
For each $\alpha<\lambda$, fix a set $U_{\alpha} \in [\lambda]^\alephes$
containing $U_h$ and a Borel function
$b_{\alpha}\colon (\NN)^{U_{\alpha}} \to P(\omega)$, coded in the ground model, such that
$$p\forces \name{A}_{\alpha} = b_{\alpha} (\seq{\name{c}_{\beta}}{\beta \in U_{\alpha}}).$$
Using the $\Delta$-System Lemma, find $W \in [\mu^{+}]^{\mu^{+}}$ and $U_{*}$ such that
$U_{\alpha} \cap U_{\beta} = U_{*}$
for all distinct $\alpha, \beta \in W$. As $U_h \sub U_{\alpha}$ for each $\alpha$, we have that $U_h \sub U_{*}$.
Fix distinct $\alpha, \beta \in W$, such that
$$\alpha \notin U_{\beta} \mbox{ and } \beta \notin U_{\alpha}.$$
This can be done as follows: Select any $\alpha \in W \setminus U_{*}$ and distinct
$\beta_0,\beta_1\in W \setminus U_{\alpha}$.
If $\alpha \notin U_{\beta_0}$, then put $\beta = \beta_0$. Otherwise,
$\alpha \notin U_{\beta_1}$ because $U_{\beta_0} \cap U_{\beta_1} = U_{*}$. In this case, put $\beta=\beta_1$.
We know that
$$p\forces \name{A}_{\alpha} \as \name{A}_{\beta} \mbox{ or } \name{A}_{\beta} \as \name{A_{\alpha}}.$$
There is $p_0 \leq p$ such that
$$(\, p_0\forces \name{A}_{\alpha} \as \name{A}_{\beta}\,)
\mbox{ or }
(\, p_0\forces \name{A}_{\beta} \as \name{A_{\alpha}} \,).$$
Without loss of generality, we may assume that
$$p_{0}\forces \name{A}_{\alpha} \as \name{A}_{\beta}.$$
Take $n_1$ and a condition $p_1 \leq p_{0} $ such that:
\be
\item $p_1 \forces \name{A}_{\alpha} \setminus n_1 \sub \name{A}_{\beta} \setminus n_1$;
\item $\alpha,\beta \in \supp (p_1)$.
\ee
Choose $n_2$ and a condition $p_2\leq p_1$ such that
\be
\item $n_2 > \max \{n_1, \max ( \dom (p_1(\beta)) \}$;
\item $p_2 \forces n_2 \in \name{A}_{\alpha}$.
\ee
We know that $ p\forces \name{A}_{\alpha} = b_{\alpha} (\smallseq{\name{c}_{\beta}}{\beta \in U_{\alpha}})$,
and thus
$\name{A}_{\alpha}$ is a $(\bbC_\lambda)_{U_{\alpha}}$-name where
$$ (\bbC_\lambda)_{U} := \set{q\in \bbC_\lambda}{\supp (q) \sub U}.$$
Thus, we may assume that
$$p_2 \rest \lambda \setminus U_{\alpha} = p_1 \rest \lambda \setminus U_{\alpha}.$$
As $\name{h}$ is a $(\bbC_\lambda)_{U_h}$-name, there are $m_{*}$ and
a condition $p_3 \leq p_2$ such that
\be
\item $p_3 \rest \lambda \setminus U_{h} = p_2 \rest\lambda \setminus U_{h}$;
\item $p_3 \forces \name{h}(n_2) =m_*$.
\ee
Finally, choose $p_4 \in \bbC_\lambda$ such that
\be
\item $\supp(p_4) = \supp (p_3)$;
\item $p_4\rest \lambda \setminus \{\beta \}= p_3 \rest \lambda \setminus \{\beta \} $;
\item $p_4(\beta) = \sigma,$ where $\sigma \colon \dom (p_{3}(\beta)) \cup \{n_2 \} \to \omega$ is defined by
$$\sigma (k)=
\begin{cases}
p_3(\beta)(k) & k\in \dom(p_3(\beta)) \\
m_{*}+1 & k = n_2.
\end{cases}
$$
\ee
In summary, the condition $p_4$ forces that
\be
\item $\name{A}_{\alpha} \setminus n_1 \sub \name{A}_{\beta} \setminus n_1 $;
\item $n_2 \in \name{A}_{\alpha}$;
\item $\name{A}_{\beta} \as \smallset{n \in \omega}{\name{c}_{\beta}(n) \leq \name{h}(n)}$;
\item $\name{h}(n_2) = m_{*}$;
\item $\name{c}_{\beta}(n_2) = m_{*}+1$.
\ee
The conditions (4) and (5) imply that
$$p_4\forces n_2 \notin \smallset{n\in \omega}{\name{c}_{\beta}(n) \leq \name{h}(n)}.$$
This is a contradiction with (2) and (3).
\epf

\begin{cor}
Let $V$ be a model of \CH{}.
For each cardinal $\lambda > \aleph_2$ of uncountable cofinality,
$$V^{\bbC_\lambda}\models\fp=\fb=\aleph_1<\lr=\lx =\aleph_2<\lambda=\cov(\cM)=\fd=\fc.$$
\end{cor}
\bpf
By Theorem \ref{thm:Cohen},
$$\lr\le\lx\le\aleph_2<\lambda=\cov(\cM)=\fd=\fc$$
in $V^{\bbC_\lambda}$.
By Theorem \ref{thm:aleph1}, $\lr>\aleph_1$ in $V^{\bbC_\lambda}$.
The remaining assertions are well known.
\epf

\subsection{A model for $\fp\ll\lr$}

Our model will be constructed using Mathias-type forcing notions.
For a centered family $\cF$ which contains all co-finite sets,
the \emph{$\cF$-Mathias forcing} is the c.c.c.\ forcing notion
$$\bbP = \set{\langle v,A \rangle\in [\omega]^{<\omega}\x \cF}{\max v< \min A}, $$
ordered by:
$$\langle u,B \rangle \leq \langle v,A \rangle \mbox{ if and only if } u \supseteq v,\ B\sub A \mbox{ and } u\setminus v\sub A.$$
This forcing notion adds a pseudointersection to the family $\cF$. Indeed, if $G$ is $\bbP$-generic, then
$$\Un\set{v}{\langle v,A \rangle\in G}$$
is a pseudointersection of $\cF$.


\bthm\label{hardmodel}
Assume the Generalized Continuum Hypothesis, and let $\mu$, $\kappa$ and $\lambda$ be uncountable cardinal numbers such that
$\kappa = \cof(\kappa) <\mu = \cof(\mu) < \lambda = \lambda^{<\mu}$.
There is a c.c.c.\ forcing notion $\bbP$ of cardinality $\lambda$ such that
$$\forces_{\bbP} \fp=\fb=\kappa < \lr =\mu< \lambda=\fc.$$
\ethm

Instead of building a model directly, as in the previous section, we will consider a transfinite sequence of
classes of forcing notions,
$\Theta_{\xi}$, and with their help we will define the forcing notion we are looking for.

A forcing notion $\bbO$ belongs to the class $\Theta_{\xi}$ if
$\bbO$ is given by an iteration $\cI$ such that:
\be
\item $\smallseq{\bbP_{\alpha}, \dQ_{\beta}}{\alpha \leq \lambda \cdot \xi, \beta < \lambda \cdot \xi}$
is a finite support iteration of length $\lambda \cdot \xi$ (ordinal product);
\item $\bbO= \bbP_{\lambda \cdot \xi}$;
\item $\bbP_0$ is the trivial forcing;
\item for each $\alpha < \lambda \cdot \xi$,
$\forces_{\bbP_{\alpha}} \dQ_{\alpha}$ is an $\name{\cF}_{\alpha}$-Mathias forcing;
\item $ \name{\cF}_{\alpha}$ is a name for a filter generated by the cofinite sets together with the family
$\smallset{\name{A}_{\alpha,\iota}}{\iota < \iota_{\alpha}}$, where
$\iota_{\alpha}$ is an ordinal $< \mu$;
\item $\iota_{\alpha}=0$ for $\alpha < \lambda$ (thus $\bbQ_{\alpha}$ is isomorphic to Cohen's forcing for $\alpha < \lambda$);
\item $\name{A}_{\alpha,\iota}$ is a $\bbP_{\alpha}$-name for a subset of $\omega$;
\item $b_{\alpha, \iota} \colon (2^{\omega})^{\omega} \to \roth $ is a Borel function from the
Cantor cube $(2^{\omega})^{\omega}$ into $\roth$, coded in the ground model;
\item $\forces_{\bbP_{\alpha}} \name{A}_{\alpha, \iota} = b_{\alpha, \iota} ( \smallseq{\name{B}_{\gamma ( \alpha, \iota, n )}}{n<\omega} )$, where  $B_{\gamma}\sub [\w]^{\w}$ denotes the $\gamma$-th generic real;
\item If $\alpha = \lambda \cdot \zeta + \nu$ (where $\nu<\lambda$), then
$$ \gamma (\alpha, \iota, n) < \lambda \cdot \zeta.$$
\item
For each $\zeta < \xi$ and each sequence $\smallseq{b_{\iota}}{\iota < \iota_*}$ of Borel functions
$b_{\iota} \colon (2^{\omega})^{\omega} \to \roth $ of length $\iota_*<\mu$,
and all ordinal numbers $\delta (\iota, n)<\lambda \cdot \zeta $
such that $\bbP$ forces that the filter generated
by the cofinite sets together with the family
$$\smallset{b_{\iota} (\smallseq{\name{B}_{\delta ( \iota, n )}}{n<\omega})}{\iota < \iota_*},$$
is proper,
there are arbitrarily large $\alpha < \lambda \cdot ( \zeta +1) $ such that:
\be
\item $\iota_{\alpha}=\iota_*$;
\item $b_{\alpha, \iota} = b_{\iota} $ for all $\iota < \iota_*$;
\item $\gamma (\alpha, \iota, n) = \delta ( \iota, n)$ for all $\iota < \iota_*$ and all $n$.
\ee
\ee

If a forcing notion $\bbO \in \Theta_{\xi}$ is obtained by an iteration
$\smallseq{\bbP_{\alpha}, \dQ_{\beta}}{\alpha \leq \lambda \cdot \xi, \beta < \lambda \cdot \xi} $, then we set
$\bbO_{\alpha} = \bbQ_{\alpha}$ for all $\alpha$.

We say that a forcing $\bbX$ is the restriction of a forcing $\bbO$ to an ordinal $\xi$,
$\bbX=\bbO\rest \xi$, if there is $\zeta \geq \xi$ such that $\bbX \in \Theta_{\xi}$, $\bbO\in \Theta_{\zeta}$ and $\bbO_{\alpha} = \bbX_{\alpha}$ for all $\alpha < \lambda \cdot \xi$.

\blem
The classes $\Theta_{\xi}$ have the following properties:
\be
\item If $\bbO \in \Theta_{\xi}$ and $\zeta <\xi$, then $\bbO\rest\zeta \in \Theta_{\zeta}$;
\item $\Theta_0$ is nonempty;
\item If $\bbO \in \Theta_{\xi}$, then there is $\bbX \in \Theta_{\xi+1}$ such that $\bbX\rest\xi =\bbO$;
\item If $\xi$ is a limit ordinal and $\smallseq{\bbO^{\zeta}}{\zeta < \xi}$ is a sequence of forcing notions
such that $\bbO^{\zeta} \in \Theta_{\zeta}$ and $\bbO^{\zeta} \rest \eta = \bbO^{\eta}$ for all $\eta < \zeta < \xi$, then
there is a
unique $\bbO^{\xi} \in \Theta_{\xi}$ such that $ \bbO^{\xi} \rest \zeta = \bbO^{\zeta}$ for all $ \zeta < \xi$.
\ee
\elem
\bpf
The only nontrivial property is (3).
To define $\bbX$, it is suffices to find functions $b_{\alpha, \iota}$
and numbers $\gamma ( \alpha, \iota, n )$ for $\alpha \in
[\lambda \cdot \xi, \lambda \cdot (\xi+1) )$ such that the conditions (10) and (11) hold.
Let $ \smallseq{\cP_{\alpha}}{\lambda \cdot \xi \leq \alpha < \lambda \cdot (\xi +1)}$ be the sequence of
all possible pairs
$$\cP =\left\langle \seq{b_{\iota}}{\iota < \iota_*}, \seq{\delta (\iota,n)}{\iota < \iota_*, n<\omega}
\right\rangle$$
where
\be
\item $\iota_*<\mu$;
\item $\smallseq{b_{\iota}}{\iota < \iota_*} $ is a sequence of Borel functions;
\item $\smallseq{\delta(\iota,n)}{\iota < \iota_*, n<\omega} $ is a matrix of ordinal
numbers $\delta (\iota,n)<\lambda \cdot \xi$;
\item the filter generated by the cofinite sets and the family
$$\smallset{
b_{\iota} (\smallseq{\name{B}_{\delta ( \iota, n )}}{n<\omega})}{\iota < \iota_*}$$
is proper.
\ee
We request that each pair appears cofinally often in this sequence.
When
$$\cP_\alpha =\left\langle
  \seq{b_{\iota}}{\iota < \iota_*}, \seq{\delta (\iota,n)}{\iota < \iota_*, n<\omega}
\right\rangle,$$
write
$b_{\alpha, \iota}= b_{\iota}$, $\gamma (\alpha, \iota, n) =\delta ( \iota, n) $ and $\iota_{\alpha}= \iota_*$.
\epf

Using the above lemma, take a sequence
$\smallseq{\bbO^{\xi}}{\xi \le \kappa}$ of forcing notions
such that $\bbO^{\xi} \in \Theta_{\xi}$ and
$\bbO^{\xi}\rest\zeta= \bbO^{\zeta}$ for every $\zeta < \xi$.
Let $\bbP_{\alpha}= \bbO^{\kappa}_{\alpha}$ for all $\alpha \le \lambda \cdot \kappa$.
The forcing notions $\bbP_{\alpha}$ are well defined:
$\bbO^{\xi}_{\alpha} = \bbO^{\kappa}_{\alpha}$ for $\xi<\ \kappa $ and $\alpha < \lambda \cdot \xi$.

\blem\label{lem:pval}
$\forces_{\bbP_{\lambda \cdot \kappa}} \fp = \fb =\kappa$.
\elem
\bpf
($\fp\geq \kappa$)
Let $\cA =\smallset{A_{\iota}}{\iota < \iota_*}\in V[G]$ be a centered family of cardinality $<\kappa$.
Let $\xi_0 <\lambda \cdot \kappa $ be such that $\cA \in V[G_{\xi_0}]$ ($\xi_0$ exists since we consider finite support iteration and $\kappa = \cof (\kappa)> \aleph_0$).
We claim that there is $\alpha > \xi_0$ such that $\cA\sub\cF_{\alpha}$.
Indeed, consider functions $b_{\iota}\colon (2^{\omega})^{\omega} \to \roth $ and ordinals $\delta(\iota, n)$ such that
$$b_{\iota} (\seq{B_{\delta( \iota, n)}}{n<\omega} ) = A_{\iota}$$
for all $\iota<\iota_*$.
By condition (11), there is $\alpha$ such that: $\iota_{\alpha}=\iota_*$,
$b_{\alpha, \iota} = b_{\iota} $ and $\delta( \iota, n) = \gamma (\alpha, \iota,n)$ for all $\iota < \iota_*$ and all $n$.
Thus, $B_{\alpha}$ is a pseudointersection of $\cA$.

($\fb\leq \kappa$) Let $f_{\xi} \in \w^{\w}$ be an enumeration of $B_{\lambda\cdot \xi}$ in $V[G]$. Then a family $\smallset{f_{\xi}}{\xi<\kappa}$ is unbounded.
\epf

\blem
$\forces_{\bbP_{\lambda \cdot \kappa}} \lr\geq \mu$.
\elem
\bpf
Assume that some
$p\in \bbP_{\lambda \cdot \xi}$ forces that a family $\cA = \smallset{A_{\iota}}{\iota < \iota_*}$, where $\iota_* < \mu$, is contained in $\roth $ and closed under finite intersections.
There are numbers $\delta (\iota, n) <\lambda \cdot \kappa$ (for $\iota < \iota_*$ and $n<\w$) and Borel functions
$b_{\iota} \colon (2^{\omega})^{\omega} \to \roth $ such that
$$p\forces_{\bbP_{\lambda \cdot \kappa}}
\name{A}_{\iota} = b_{\iota}(\seq{\name{B}_{\delta(\iota, n)}}{n<\omega} ) \mbox{ for all } \iota.$$
Write each $\delta(\iota,n)$ in the form $$\delta(\iota,n)= \lambda \cdot \zeta(\iota,n) + \eta(\iota,n),$$
where $\eta(\iota,n) < \lambda$.
Set $$\eta_* = \sup \set{\eta(\iota,n)}{\iota<\iota_*, n<\omega}.$$
As $\cof(\lambda)\ge\mu$ (indeed, $\lambda = \lambda^{<\mu}$),
we have that $\eta_*< \lambda$.
Since, in addition, $\mu$ is regular,
there is $S \sub\lambda \cdot \kappa$ of cardinality $<\mu$ such that
\be
\item $\smallset{\delta (\iota,n)}{\iota < \iota_{*},\ n <\omega}\sub S$;
\item If $\alpha \in S$ then
$\smallset{\gamma (\alpha, \iota,n)}{\iota < \iota_{\alpha},\ n <\omega} \sub S$.
\ee
Set
$$ S_{\xi} = S \cap \lambda\cdot \xi.$$
Then
$\smallseq{S_{\xi}}{\xi < \kappa}$ is a $\sub $-increasing sequence.
Let $U_\xi=\smallset{\iota<\iota_*}{\forall n,\ \delta(\iota,n)\in S_\xi}$, so that $\smallseq{U_\xi}{\xi<\kappa}$ is
$\sub$-increasing with union $\iota_*$.

Choose $\beta_{\xi}$ and $\eta_{\xi}$, $\xi < \kappa$, by induction such that
\be
\item $\beta_{\xi} = \lambda \cdot \xi + \eta_{\xi}$ where $\eta_{\xi} <\lambda$;
\item $\beta_{\xi} \notin S$;
\item $\eta_{\xi} > \sup ( \smallset{\eta_{\zeta}}{\zeta < \xi} \cup \{\eta_* \} ) $;
\item $\forces_{\bbP_{\lambda \cdot \xi}} \smallset{\name{A}_{\beta_\xi,\iota}}{\iota<\iota_{\beta_\xi}} =
\smallset{\name{A}_{\iota}}{\iota\in U_\xi} \cup \smallset{\name{B}_{\beta_{\zeta}}}{\zeta < \xi}.$
\ee
The induction can be carried out,  since $$\forces_{\bbP_{\lambda \cdot \kappa}}\name{B}_{\beta_{\xi}} \as \name{B}_{\beta_{\zeta}} \mbox{ for } \zeta < \xi $$
and
$$\forces_{\bbP_{\lambda \cdot \kappa}} \name{B}_{\beta_{\xi}} \mbox{ has infinite intersection with every member of } \cA. $$
To observe that the last condition holds, it suffices to use (4) and the fact that $\beta_{\xi} \notin S$.

Since $\forces_{\bbP_{\lambda \cdot \kappa}} \bigcup_{\xi<\kappa} U_{\xi} =\iota_{*}$,
we conclude by (4) and the definition of $\name{\bbQ}_{\beta_{\xi}}$ that,
in $V^{\bbP_{\lambda\cdot \kappa}}$,
the set $\set{B_{\beta_{\xi}}}{\xi<\kappa}$ is a linear refinement of
$\set{A_{\iota}}{\iota<\iota_{*}}$.
\epf

\blem\label{lem:nopi}
If $U \in [\lambda]^{\mu}$ is from the ground model $V$ and $\gamma \in [\lambda, \lambda \cdot \kappa ],$ then
$$(\star) \ \forces_{\bbP_{\gamma}} \bigl| \smallset{\alpha \in U}{\name{C} \as \name{B}_{\alpha}}\bigr| < \mu \mbox{ for each infinite } \name{C} \sub \omega.$$
\elem
\bpf
We prove the fact by induction on $\gamma \in [\lambda, \lambda \cdot \kappa ]$. For each $\gamma$ let $G_{\gamma}$ denote the $\bbP_{\gamma}$-generic filter.

Assume that  $\gamma = \lambda$. Let $C\in V[G_{\lambda}]$ be an infinite subset of $\omega$ and let $\name{C}$ be a $\bbP_{\lambda}$-name for $C$. As $C$ is determined by countably many Cohen reals, we may assume by changing
the order that $C\in V[G_\omega]$.
Then
$$\forces_{\bbP_{\gamma}}
\bigl| \smallset{\alpha \in U}{\name{C} \as \name{B}_{\alpha}}\bigr| < \aleph_1.$$
We next establish the preservation of the condition $(\star)$ through the steps of iteration.

Assume that $\gamma = \beta +1$ is a successor ordinal.
We will work in $V[G_{\beta}]$ and force with $\bbQ_{\beta}$.
Assume that in $V[G_{\beta}]$, for every infinite $C$
  $$ (\star) \  |\set{\alpha\in U}{C\as B_{\alpha}}|<\mu.$$
We force with $\cF_{\beta}$-Mathias forcing $\bbQ_{\beta}$, where  $\cF_{\beta}$ is generated by a centered family of cardinality $<\mu$.
Therefore $\bbQ_{\beta}$ contains a dense subset $D$ of cardinality $<\mu$.
Assume that
$$\forces_{\bbQ_{\beta}}  |\set{\alpha\in U}{\name{C}\as B_{\alpha}}| \geq\mu \mbox{ for some infinite } \name{C}\sub \w.$$
Let
$$ W=\{ \alpha\in U: (\exists q_{\alpha} )  \ q_{\alpha}\forces \name{C}\as B_{\alpha} \}.$$
The set $W$ belongs to the ground model $V[G_{\beta}]$.
We may assume that for each $\alpha\in W$, $q_{\alpha} \in D$ .
By pigeonhole principle there is $q_{*} \in D $ and a set $W_1\sub W$ of cardinality $\mu$ such that
$q_{\alpha}= q_*$ for each $\alpha \in W_1$.
This means that for each $\alpha \in W_1$
$$q_* \forces \name{C} \as B_{\alpha} .$$
For each $\alpha \in W_1$ there is $r_{\alpha} \leq q_*$ and $k_{\alpha}$ such that for each  $\alpha \in W_1$:
$$r_{\alpha} \forces \name{C} \setminus [0,k_{\alpha}) \sub B_{\alpha}.$$
Again, by pigeonhole principles there are $r_*$ and $k_*$ and $W_2 \sub W_1$ of cardinality $\mu$ such that
$r_{\alpha}= r_*$, $k_{\alpha}=k_*$ for each $\alpha \in W_2$.
This means that  for each $\alpha \in W_2$:
$$r_{*} \forces \name{C} \setminus [0,k_{*}) \sub B_{\alpha}.$$
It follows that
$$r_{*} \forces \bigcap_{\alpha \in W_2} B_{\alpha} \mbox{ is infinite  }.$$
But $\bigcap_{\alpha \in W_2} B_{\alpha}$ belongs to $V[G_{\beta}]$ contradicting $(\star)$.

Assume that  $\gamma$ is a limit ordinal of uncountable cofinality.
Let $C\in V[G_{\lambda}]$ be an infinite subset of $\omega$ and let $\name{C}$ be a $\bbP_{\gamma}$-name for $C$.
By Lemma 16.14 in \cite{Je} there is $\beta < \gamma$ such that $C\in V[G_{\beta}]$ and $\beta \geq \lambda$.
By the inductive hypothesis,  we have that
$$\forces_{\bbP_{\beta}} \bigl| \smallset{\alpha \in U}{\name{C} \as \name{B}_{\alpha}}\bigr| < \mu.$$
Since c.c.c.\ forcing notions preserve cardinality, we have that
$$\forces_{\bbP_{\gamma}} \bigl| \smallset{\alpha \in U}{\name{C} \as \name{B}_{\alpha}}\bigr| < \mu.$$
Finally, let  $\gamma>\lambda$ be a limit ordinal with countable cofinality. Fix a sequence $\smallseq{\gamma_n}{n<\w}$ increasing to $\gamma$.
Towards a contradiction assume that there is $p\in \bbP_{\gamma}$ such that $$p \forces_{\bbP_{\gamma}} \name{\cU} = \smallset{\alpha \in U}{\name{C} \as \name{B}_{\alpha}} \mbox{ has cardinality }\mu.$$
Let $\name{\beta}_{\iota}$ be a name for the $\iota$-th element of $\name{\cU}$, so that
Let $G$ be a $\bbP_{\gamma}$-generic filter  containing $p$ and for every $\iota <\mu$ let $p_{\iota} \in \bbP_{\gamma}$, $\alpha_{\iota}<\mu$ and $k_{\iota}<\w$ be such that
\begin{enumerate}
\item $p_{\iota} \leq p$;
\item $p_{\iota} \forces \name{\beta}_{\iota} = \alpha_{\iota} $;
\item $p_{\iota} \forces \name{C} \setminus \name{B}_{\alpha_{\iota}} \sub [0, k_{\iota})$.
\end{enumerate}
As  $\supp (p_{\iota})$ is finite for each $\iota<\mu$, there exists $n_{\iota} <\w$ such that $\supp(p_{\iota}) \sub \gamma_{n_{\iota}}$.
Since there are $\mu$ many indices $\iota$ and only countably many $n_{\iota}$ and $k_{\iota}$, there exist $n_{*}$ and $k_{*}$  such that the set
$$ W=\set{\iota <\mu}{n_{\iota} = n_{*}, k_{\iota} = k_{*}}$$
has cardinality $\mu$.
In particular we have $p_{\iota} \in G_{\gamma_{n_*}}$ for all $\iota\in W$. Notice that $W, \seq{p_{\iota}}{\iota \in W} \in V[G]$  (in fact they belong to $V[G_{\gamma_{n_*}}]$)
Let $\name{D}$ be a $\bbP_{\gamma_{n_*}}$ -name defined as follows: given $\bbP_{\gamma_{n_*}}$-filter $H$, let $\name{D}[H]$ be a set
$$\set{k\in \w}{ (\exists q\in \bbP_{\gamma_{n_*}})\ q\forces k\in \name{C}' }, $$
$\name{C}'$ is the $\bbP_{\gamma_{n_*}, \gamma}$-name obtained in a standard way by "partially evaluating $\name{C}$ with $H$".
We claim that  $p_{\iota} \forces \name{D}\setminus \name{B}_{\alpha_{\iota}} \sub [0,k_*)$ for all $\iota$.
Indeed, otherwise there exists $r\leq p_i$,$r\in \bbP_{\gamma_{n_*}}$ and $k>k_*$ such that $r\forces_{\bbP_{\gamma_{n_*}}} k\in \name{D} \setminus \name{B}_{\alpha_{\iota}}$. Thus
$$r\forces_{\bbP_{\gamma_{n_*}}} \bigl( (  \exists q\in \bbP_{\gamma_{n_*},\gamma}) \ q\forces_{\bbP_{\gamma_{n_*},\gamma}} k\in \name{C}' \setminus \name{B}_{\alpha_{\iota}}   \bigr).$$
Let $r' \leq r$ $r' \in  \bbP_{\gamma_{n_*}}$ and $q\in \bbP_{\gamma_{n_*},\gamma}$ be such that
$$r\forces_{\bbP_{\gamma_{n_*}}} \bigl(    q\forces_{\bbP_{\gamma_{n_*},\gamma}} k\in \name{C}' \setminus \name{B}_{\alpha_{\iota}}   \bigr).$$
This means that $r'^{\frown} q \forces_{\bbP_{\gamma}} k\in \name{C}\setminus \name{B}_{\alpha_{\iota}} $. But this is impossible because
$$r'^{\frown} q \leq p_{\iota} \forces_{\bbP_{\gamma}} \name{C}\setminus \name{B}_{\alpha_{\iota}} \sub k_{\iota}.$$
This proves the claim.

As $p_\iota\in G_{\gamma_{n_*}}$, the above claim implies that, in $V[G]$, 
$$W\sub \smallset{\iota}{\name{D}[G_{\gamma{n_*}}] \sub \name{B}_{\alpha_{\iota}}[G_{\gamma{n_*}}]}.$$
Take $p'\leq p$,$p'\in \bbP_{\gamma_{n_*}}$ that forces this inclusion.
Then $$p'\forces \bigl| \smallset{\alpha }{\name{D} \sub B_{\alpha} }\bigr|=\mu,$$ contradicting the inductive hypothesis.
\epf

\blem\label{lemacik}
$\forces_{\bbP_{\lambda \cdot \kappa}} \lr \le \mu$.
\elem
\bpf
Assume otherwise. Then there is $p\in \bbP$ such that $p\forces \mu <\lr$. Take a generic filter $G$ containing $p$.
We argue in $V[G]$.
The family $\cB= \smallset{B_{\alpha}}{\alpha < \mu}$ has a linear refinement. By Lemma 2.4, either there is a tower $\smallset{T_{\iota}}{\iota< \fp}$ refining $\cB$, or $\cB$ has a pseudointersection.
The second case cannot happen since it contradicts $(\star)$ of Lemma \ref{lem:nopi}.
Let $D_{\iota}= \smallset{\alpha}{T_{\iota}\as B_{\alpha}}$.
Since $\bigcup_{\iota <\fp} D_{\iota} = \mu$, there is $\iota$ such that $\card{D_{\iota}}=\mu$.
This means that $T_{\iota}\as B_{\alpha}$ for $\mu$-many $\alpha$, contradicting $(\star)$ of Lemma \ref{lem:nopi}.
\epf

This completes the proof of Theorem \ref{hardmodel}. \qed

\subsection{A model for $\lr \ll \fb=\lx=\fd \ll \fc$}
Our model will be constructed using Mathias-type forcing notions as in the previous section,
together with Hechler forcing.

\bthm\label{thm:Hechler}
Assume the Generalized Continuum Hypothesis, and let $\kappa$, $\eta$ and $\lambda$ be uncountable cardinal numbers such that
$\kappa$ and $\eta$ are regular and $\kappa < \eta < \lambda=\lambda^{<\kappa}$.
There is a c.c.c.\ forcing notion $\bbP$ of cardinality $\lambda$ such that
$$\forces_\bbP \fp =\kappa <  \lr = \kappa^{+}  \leq \fb = \lx = \fd =\eta < \lambda=\fc.$$
\ethm
\bpf
We use the iteration of Theorem \ref{hardmodel}, but $\kappa$, $\eta$ and $\lambda$ here stand for
$\mu$, $\kappa$, and $\lambda$ there, respectively,
and we intersperse Hechler's forcing during the iteration. More precisely, the forcing notion $\bbP$ is
given by the following iteration:
\be
\item $\smallseq{\bbP_{\alpha}, \dQ_{\beta}}{\alpha \leq \lambda \cdot \eta, \beta < \lambda \cdot \eta}$
is a finite support iteration of length $\lambda \cdot \eta$ (ordinal product);
\item $\bbP= \bbP_{\lambda \cdot \eta}$;
\item $\bbP_0$ is the trivial forcing;
\item If  $\alpha  \in \smallset{\lambda \cdot \xi}{ \xi >0}$, then $\forces_{\bbP_{\alpha}} \dQ_{\alpha}$ is Hechler's forcing;
\item If $\alpha<\lambda\cdot\eta$ and $\alpha  \notin \smallset{\lambda \cdot \xi}{ \xi >0}$, then
\be
\item $\forces_{\bbP_{\alpha}} \dQ_{\alpha}$ is an $\name{\cF}_{\alpha}$-Mathias forcing;
\item $ \name{\cF}_{\alpha}$ is a name for a filter generated by centered family
$\smallset{\name{A}_{\alpha,\iota}}{\iota < \iota_{\alpha}}$ which contains cofinite sets, where
$\iota_{\alpha}$ is an ordinal $< \kappa$;
\item $\iota_{\alpha}=0$ for $\alpha < \lambda$ (thus $\bbQ_{\alpha}$ is isomorphic to Cohen's forcing for $\alpha < \lambda$);
\item $\name{A}_{\alpha,\iota}$ is a $\bbP_{\alpha}$-name for a subset of $\omega$;
\item $b_{\alpha, \iota} \colon (2^{\omega})^{\omega} \to \roth $ is a Borel function coded in the ground model;
\item $\forces_{\bbP_{\alpha}} \name{A}_{\alpha, \iota} = b_{\alpha, \iota} ( \smallseq{\name{B}_{\gamma ( \alpha, \iota, n )}}{n<\omega} )$, where $B_{\alpha}\sub [\w]^{\w}$ denotes the $\alpha$-th generic real;
\item If $\alpha = \lambda \cdot \xi + \nu$ (where $\nu<\lambda$), then
$$ \gamma (\alpha, \iota, n) < \lambda \cdot \xi.$$
\item
For each $\zeta < \xi$ and each sequence $\smallseq{b_{\iota}}{\iota < \iota_*}$ of Borel functions
$b_{\iota} \colon (2^{\omega})^{\omega} \to \roth $ of length $\iota_*<\kappa$,
and all ordinal numbers $\delta (\iota, n)<\lambda \cdot \zeta $
such that $\bbP$ forces that the filter generated
by the cofinite sets together with the family
$$\smallset{b_{\iota} (\smallseq{\name{B}_{\delta ( \iota, n )}}{n<\omega})}{\iota < \iota_*},$$
is proper,
there are arbitrarily large $\alpha < \lambda \cdot ( \zeta +1) $ such that:
\be
\item $\iota_{\alpha}=\iota_*$;
\item $b_{\alpha, \iota} = b_{\iota} $ for all $\iota < \iota_*$;
\item $\gamma (\alpha, \iota, n) = \delta ( \iota, n)$ for all $\iota < \iota_*$ and all $n$.
\ee
\ee
\ee

Observe that $\forces_{\bbP} \fb = \lx = \fd =\eta$ since Hechler reals are  added in steps $\lambda\cdot\xi$ ($\xi<\eta$) of the iteration. Also, $\forces_{\bbP} 2^\alephes=\lambda$ holds, since $\lambda=\lambda^{<\kappa}$.
It remains to prove
that $\forces_{\bbP} \fp = \kappa$  and  $\forces_{\bbP} \lr = \kappa^{+}$.

\blem\label{lemcomb2}
If $U \in [\lambda]^{\kappa}$ is from the ground model $V$ and $\gamma \in [\lambda, \lambda \cdot \kappa ],$ then
$$(\star) \ \forces_{\bbP_{\gamma}} \bigl| \smallset{\alpha \in U}{\name{C} \as \name{B}_{\alpha}}\bigr| < \kappa \mbox{ for each infinite } \name{C} \sub \omega.$$
\elem

\bpf
The proof is as in Lemma \ref{lem:nopi}, with one more case to check:
$\gamma=\beta+1$ and $\forces_{\bbP_{\beta}} \name{\bbQ}$
is Hechler's forcing.

In $V$, enumerate $U =\smallset{ \alpha_{\delta}}{\delta <\kappa}$.
Consider a family
$\smallset{B_{\alpha}}{\alpha \in U} =
\smallset{B_{\alpha_{\delta}}}{\delta<\kappa}$ in $V[G_{\beta}]$.
It is \emph{eventually narrow}, that is, for each $C\in \roth$ there is $\delta_0$ such that $C\nsubseteq^* B_{\alpha_{\delta}}$ 
for each $\delta >\delta_0$.
By \cite[Theorem 3.1]{BaumDord}, eventually narrow families are preserved
by Hechler's forcing.
Thus,
$$\forces_{\bbP_{\gamma}} \bigl| \smallset{\alpha \in U}{\name{C} \as \name{B}_{\alpha}}\bigr| < \kappa \mbox{ for each infinite } \name{C} \sub \omega.\qedhere$$
\epf

\blem
$\forces_{\bbP} \kappa \leq \fp$.
\elem
\bpf
The proof is as in Lemma \ref{lem:pval}.
\epf

\blem
$\forces_{\bbP} \fp \leq \kappa$.
\elem
\bpf
By Lemma \ref{lemcomb2} for $U=\kappa$,
the family $\smallset{B_{\alpha}}{\alpha \leq \kappa}$ of first $\kappa$
Cohen reals is an example of a centered family in $V[G]$
that has no pseudointersection.
\epf

\blem
$\forces_{\bbP} \lr \le \kappa^{+}$.
\elem

\bpf
The proof is as in Lemma 3.8. The only difference is that since now $\kappa =\mu= \fp$,
we need to change $\kappa$ to $\kappa^+$ in the conclusion.
\epf

Consider the following weak version of the Martin's Axiom $M(\kappa)$: If \begin{enumerate}
\item $\cA \sub \roth$ is a centered family of cardinality $<\kappa$ (where $\kappa >\w$), that contains all cofinite sets,
\item $\bbQ=\bbQ_\cA$ is the $\cA$-Mathias forcing notion,
\item $\cD_{\beta}$ is an open dense subset of $\bbQ$ for each $\beta <\kappa$;
\end{enumerate}
 then there is a filter $H\sub \bbQ$ such that $H\cap \cD_{\beta} \neq \emptyset$ for each $\beta < \kappa$.

\blem
$M(\kappa)$ implies that $\lr \geq \kappa^{+}$.
\elem

\bpf
Assume that $\smallset{ A_{\alpha} }{\alpha < \kappa }$ is a centered family. We may assume  that it contains all cofinite set and is closed under finite intersection.
We choose $A^{-}_{\alpha}$ by induction on $\alpha <\kappa$ such that:
\begin{enumerate}
\item $A^{-}_{\alpha} \as A^{-}_{\beta}$ for each $\beta <\alpha$,
\item $ A^{-}_{\alpha} \cap A_{\beta}$ is infinite for each $\beta<\kappa$.
\end{enumerate}
Assume that $A^{-}_{\beta}$ is defined for $\beta < \alpha$. Let $\cA$ be the closure of the family $\smallset{A^{-}_{\beta}}{\beta <\alpha}$ under finite intersections and cofinite sets. Apply $M(\kappa)$ to the
family $\cA$ and the dense sets
$\cD_{\beta,k} = \smallset{ \langle u, B\rangle }{ \card{ u \cap A_{\beta} } \geq k  }$ (where $\beta <\kappa$ and $n<\w$) to obtain $H$.
Set $A^{-}_{\alpha}=\Un \smallset{u}{\la u,B\ra\in H}$.
\epf

\blem
$\forces_{\bbP} \lr \geq \kappa^{+}$.
\elem

\bpf
Assume that $\cA$ is a centered family of cardinality $<\kappa$ in the extended model $V[G]$,
which contains all cofinite sets,
and $\mathscr{D}$ is a family of $\le\kappa$
open dense subsets of $\bbQ=\bbQ_\cA$.
Assume that $p$ forces that $\name{\cA}=\smallset{\name{A}_\iota}{\iota<\name{\iota_*}<\kappa}$ and
$\name{\mathscr{D}}=\smallset{\name{\cD}_{\epsilon}}{\epsilon<\kappa}$
form a counter-example.
The forcing is c.c.c., and $p$ forces that $\name{\iota_*}<\kappa$.
We may assume that $\iota_*$ is in the ground model.

As $\eta=\cof(\eta)>\kappa$, we may assume that all $\name{A}_\iota, \name{\cD}_{\epsilon}$
are $\bbP_{\lambda\cdot\xi}$-name for some $\xi<\eta$.
We can find $\alpha\in[\lambda\cdot\xi,\lambda\cdot(\xi+1))$ such that
$\smallset{\name{A}_{\alpha,\iota}}{\iota<\iota_\alpha}=\smallset{\name{A}_\iota}{\iota<\iota_*}$
is forced. We conclude as in the proof of the consistency of Martin's Axiom.
\epf

The proof of Theorem \ref{thm:Hechler} is completed.
\epf

\section{Open problems}

One of our main results (Theorem \ref{thm:cflr}) is that the cofinality of
$\lr$ is uncountable.

\bprb
Is it consistent that $\lr$ is singular?
\eprb

We introduce below two ad-hoc names for combinatorial cardinal characteristics. Once progress is made on the associated problems, better names may be introduced.

\bdfn
Let $\kappa_1$ be the minimal cardinality of a family $\cA\sub (\roth)^\N$ such that:
\be
\item For each $n$, $\smallset{A(n)}{A\in\cA}$ is linear.
\item There is \emph{no} $g\in\NN$ such that
the sets $S_A:=\smallset{n}{g(n)\in A(n)}$ are infinite, and the family
$\smallset{S_A}{A\in\cA}$ has a linear refinement.
\ee
\edfn

The following assertions are proved exactly as in Section \ref{sec:ZFCSPM}.

\blem
\mbox{}
\be
\item $\non(\sone(\Tau^*,\Tau^*))=\kappa_1$.
\item $\min\{\fb,\fs,\cov(\cM)\}\le \kappa_1$.
\item $\min\{\max\{\lr,\min\{\fb,\fs\}\},\cov(\cM)\}\le \kappa_1$. \qed
\ee
\elem

\bprb
Can we express $\kappa_1$, in ZFC, in terms of classic combinatorial cardinal characteristics of the continuum?
\eprb

\bdfn
Let $\kappa_{\mathrm{fin}}$ be the minimal cardinality of a family $\cA\sub (\roth)^\N$ such that:
\be
\item For each $n$, $\smallset{A(n)}{A\in\cA}$ is linear.
\item There are \emph{no} finite sets $F_0,F_1,\dots\sub\w$ such that
the sets $S_A:=\Un_n \{n\}\x (A(n)\cap F_n)\sub\w\x\w$ are infinite, and the family
$\smallset{S_A}{A\in\cA}$ has a linear refinement.
\ee
\edfn

We have the following.
\blem
\mbox{}
\be
\item $\non(\sfin(\Tau^*,\Tau^*))=\non(\sfin(\Tau,\Tau^*))=\kappa_{\mathrm{fin}}$.
\item $\min\{\cov(\cM),\lr\},\min\{\fb,\fs\}\le\kappa_{\mathrm{fin}}\le\lx$.\qed
\ee
\elem

\bprb
Can we express $\kappa_{\mathrm{fin}}$, in ZFC, in terms of classic combinatorial cardinal characteristics of the continuum?
If not, can we improve the above bounds in ZFC?
\eprb

\ed
\end{document}